\documentclass{amsart}  
\usepackage[latin1]{inputenc}
\usepackage[T1]{fontenc}
\usepackage{pstricks,geometry,amsmath,amssymb,epsfig,graphicx,graphics,float,multicol,fix-cm}

\newtheorem{theorem}{Theorem}[section]
\newtheorem{proposition}[theorem]{Proposition}
\newtheorem{lemma}[theorem]{Lemma}
\newtheorem{corollary}[theorem]{Corollary}

\title{Competitive exclusion for chemostat equations with variable yields} 

\author{Tewfik Sari} 
\address{Tewfik Sari,  
Universit\'e de Haute Alsace, LMIA,
4 rue des fr\`eres Lumi\`ere, 68093 Mulhouse, 
\&
Modemic, Inra-Inria, UMR Mistea, 2 place Viala, 34060 Montpellier, France.}
              \email{tewfik.sari@uha.fr}

\date{\today}

\begin{document}

\begin{abstract}
In this paper, we 
study the global dynamics of
a chemostat model with a single nutrient and several competing species.
Growth rates are not required to be proportional to food uptakes.
The model was studied by Fiedler and Hsu [J. Math. Biol. (2009) 59:233-253].
These authors prove 
the nonexistence of periodic orbits,
by means of a multi-dimensional Bendixon-Dulac criterion.
Our approach is based on the construction of Lyapunov functions. 
The Lyapunov functions 
extend those used by 
Hsu [SIAM J. Appl. Math. (1978) 34:760-763] and  by
Wolkowicz and Lu [SIAM J. Appl. Math. (1997) 57:1019-1043]
in the case when growth rates are proportional to food uptakes.\\

\end{abstract}
\maketitle
\section{Introduction and main result}
The aim of this paper is to prove the Competitive Exclusion Principle (CEP) for the following 
competition for a single resource model
\begin{equation}\label{eqfpi}
\begin{array}{lcll} 
 \displaystyle S'& =&D(S^0-S)-\displaystyle
\sum_{i=1}^Np_i(S)x_i&\\[3mm]
 \displaystyle x'_i&=& f_i(S)x_i,& i=1\cdots N,
\end{array}
\end{equation}
where $S(t)$ and $x_i(t)$, $i=1\cdots N$, denote respectively 
the nutrient concentration and the concentration of the $i$th 
competing species at time $t$. 
The input concentration $S^0$ and the removal rate $D$ are assumed to be constant. 
The uptake rate $p_i(S)$ satisfies 
\begin{equation}\label{pi}
p_i(0)=0\mbox{ and }p_i(S)>0\mbox{ for } S>0.
\end{equation}
The growth rate $f_i(S)$ satisfies 
\begin{equation}\label{fi}
f_i(S)<0\mbox{ for }0\leq S<\lambda_i\mbox{ and }f_i(\lambda_i)=0.
\end{equation}
The smallest positive value of the concentration substrate $S=\lambda_i$ given in (\ref{fi}), 
where the growth rate $f_i(S)$ of $x_i$ is 0, is called the {\em break-even concentration} for the $i$th species. 
This model was considered by Fiedler and Hsu \cite{FH} as an extension of the classical chemostat model
\begin{equation}\label{eqsxi}
\begin{array}{lcll} 
 \displaystyle S'& =&D[S^0-S]-\displaystyle
\sum_{i=1}^N\frac{q_i(S)}{Y_i}x_i,&\\[3mm]
 \displaystyle x'_i&=& [q_i(S) - D_i]x_i,& i=1\cdots N,
\end{array}
\end{equation}
where $Y_i$ is the conversion constant, or {\em yield}, 
for the $i$th species, 
and the growth rate $q_i(S)$ satisfies 
\begin{equation}\label{qi}
q_i(0)=0\mbox{ and }q_i(S)>0\mbox{ for }S>0.
\end{equation} 
Thus, the break-even concentration $S=\lambda_i$ satisfies the condition $q_i(S)=D_i$;
it is the smallest value, where the growth $q_i(S)$ of $x_i$ is balanced by the removal rate $D_i$. 
The chemostat occupies a central place in mathematical ecology, see the monograph of Smith and Waltman \cite{chem}. 
It is a model of the dynamics of microbial competition. Basically, the chemostat consists of a 
nutrient input, 
pumped at a constant rate into a well-mixed culture vessel.
The culture vessel contains the microorganisms that are growing and competing for the nutrient.
Volume is kept constant by pumping the mixed contents out at the same rate.
If the growth functions $f_i(S)$ 
and the uptake functions $p_i(S)$ are of the form
\begin{equation}\label{proportionality}
f_i(S)=q_i(S)-D_i,\qquad p_i(S)=q_i(S)/Y_i,
\end{equation}
then the model (\ref{eqfpi}) of Fiedler and Hsu \cite{FH} reduces to the simple chemostat model (\ref{eqsxi}).
However, the particular form (\ref{proportionality}) of the growth function $f_i(S)$ is not assumed 
in \cite{FH} and hence, the model (\ref{eqfpi}) recovers the more general case of {\em variable yields} model
\begin{equation}\label{eqsqi}
\begin{array}{lcll} 
 \displaystyle S'& =&D[S^0-S]-\displaystyle
\sum_{i=1}^Np_i(S)x_i,&\\[3mm]
 \displaystyle x'_i&=& [q_i(S) - D_i]x_i,& i=1\cdots N,
\end{array}
\end{equation}
where uptake rate $p_i(S)$ satisfies (\ref{pi}) and growth rate $q_i(S)$ satisfies (\ref{qi}).
The functions $y_i(S)$, defined by 
\begin{equation}\label{yields}
y_i(S)=\frac{q_i(S)}{p_i(S)},\qquad i=1\cdots N,
\end{equation}
are the growth yields. The model (\ref{eqsqi}) was
considered by Arino, Pilyugin and Wolkowicz \cite{apw} as a generalization 
of the constant yields case (\ref{eqsxi}), for which  the yields (\ref{yields}) satisfy $y_i(S)=Y_i$. It was further investigated by Sari and Mazenc \cite{sari1,sari2}.
Without loss of generality (see Section \ref{PR}), we assume that $D=1$ and $S^0=1$ in (\ref{eqfpi}). The system becomes
\begin{equation}\label{D=1}
\begin{array}{lcll} 
 \displaystyle S'& =&1-S-\displaystyle
\sum_{i=1}^Np_i(S)x_i,&\\[3mm]
 \displaystyle x'_i&=& f_i(S)x_i,& i=1\cdots N.
\end{array}
\end{equation}

Coexistence of the $N$ species is a fundamental question on the model (\ref{D=1}) 
of competition for a single resource. 
Looking for coexistence at positive equilibria we have to solve equations $f_i(S)=0$ 
simultaneously for all $i=1\cdots N$.
In general, for $N\geq 2$, these equations cannot be solved for the same value of $S$. 
Thus, generically, (\ref{D=1}) can have the following equilibria: the washout equilibrium
\begin{equation}\label{washout}
E_0=(1,0,\cdots,0),
\end{equation}
where all species go extinct, and equilibria 
$E_{i}$, $i=1\cdots N$,
where all components of $E_{i}$ vanish, except for the first and the $(i+1)$th, which are
$$S=S^*,\qquad x_i=\frac{1-S^*}{p_i(S^*)},$$ 
where $S^*\in]0,1[$ satisfies $f_i(S^*)=0$. 
Hence, at any equilibrium point $E_{i}$, all but one species go extinct. 

Since $f(\lambda_i)=0$, 
the break-even concentration $S^*=\lambda_i$ gives rise to an equilibrium point $E_i$ for the system, if and only if $\lambda_i<1$. 
A well-known open-problem in the theory of the chemostat is to prove the global asymptotic 
stability of the equilibrium point $E_{i}$ with the lowest break-even concentration. 
If this equilibrium is globally asymptotically stable (GAS), then the CEP holds: 
only one species survives, namely the species which makes optimal use of the resource. 
The reader is referred to \cite{tilman}, for complements
and details on the CEP, and to \cite {MGPM} for recent results and a discussion on 
competitive exclusion. 
Most of the results on the CEP for (\ref{eqsxi}) and (\ref{eqsqi}) have been based on Lyapunov functions \cite{amg,hsu,li,sari1,sari2,wl,wx} .
For a survey of constructing Lyapunov functions in the chemostat, the reader is referred
to \cite{hsu1}. We simply recall here that Hsu \cite{hsu} 
proved the CEP for the Monod case of (\ref{eqsxi}), when the growth functions are
\begin{equation}\label{monod}
q_i(S)=\frac{a_iS}{b_i+S},
\end{equation}
and Wolkowicz and Lu \cite{wl} extended the result of \cite{hsu} to more general growth functions.

Instead of a Lyapunov function approach, Fiedler and Hsu \cite{FH}
applied a multi-dimensional Bendixon-Dulac criterion to exclude periodic solutions. 
Under some technical conditions on 
the functions $f_i$ and $p_i$ they proved that (\ref{D=1}) does not possess positive non-stationary periodic orbits.
In our previous works \cite{sari1,sari2}, we showed that
both Lyapunov functions used by Hsu \cite{hsu} and Wolkowicz and Lu \cite{wl} 
can be extended to the variable yields case model (\ref{eqsqi}). 
The aim of this paper is to show that these Lyapunov functions can 
also be used to obtain the CEP for (\ref{D=1}).

A necessary condition to avoid washout 
of the species, and global convergence towards the washout equilibrium $E_0$ defined by (\ref{washout}), 
is that $\lambda_i<1$ for at least one species.
Assume that the species are labeled so that $0<\lambda_1<1$. Then 
\begin{equation}\label{equilibrium}
E_{1}^*=(\lambda_1,x_{1}^*,0,\cdots,0),
\end{equation}
where $x_1=x_{1}^*=P_1(\lambda_1)$ is an equilibrium. Here
\begin{equation}\label{P1}
P_1(S)=\frac{1-S}{p_1(S)}.
\end{equation}
Using linearization of (\ref{D=1}) about $E_1^*$ one proves that:
\begin{lemma}\label{LAS}
The equilibrium (\ref{equilibrium}) is locally exponentially stable if and only if $f'_1(\lambda_i)>0$ and $P_1'(\lambda_1)<0$.
\end{lemma}
We consider the global asymptotic stability of $E_{1}^*$.
Our main result is
\begin{theorem}\label{ourthm1}
Assume that $\lambda_1<1$ and for all $0<S<1$,
\begin{equation}\label{h11}
(S-\lambda_1)f_1(S)>0,\mbox{ for }S\neq\lambda_1,
\end{equation}
\begin{equation}\label{h31}
(S-\lambda_1)(P_1(S)-P_1(\lambda_1))<0,\mbox{ for }S\neq\lambda_1.
\end{equation}
Assume that there exist constants
$\alpha_i>0$ for each $i\geq 2$ satisfying $\lambda_i<1$, such that for all $0<S<1$,
\begin{equation}\label{h21}
f_1(S)p_i(S)>\alpha_if_i(S)(1-S).
\end{equation} 
Then the equilibrium $E_{1}^*$ is GAS
for (\ref{D=1}) with respect to the interior of the positive cone.
\end{theorem}
Notice that the following property holds.
\begin{lemma}\label{lemmaCEP}
The conditions $\lambda_1<1$ and (\ref{h21}) imply that $\lambda_1<\lambda_i$ for all $i\geq 2$.
\end{lemma}
{\em Proof}
Assume that there exists $i\geq 2$ such that $\lambda_i<\lambda_1$. Then, there exists $S\leq\lambda_1$ such that $f_i(S)>0$. 
Hence, since $S\leq\lambda_1<1$, $f_i(S)(1-S)>0$. On the other hand,
using (\ref{fi}), $f_1(S)\leq 0$. Hence, the inequality (\ref{h21}) is violated. 
\medskip

\noindent
This lemma shows that the winning species $x_1$ of Theorem \ref{ourthm1} has the lowest break-even concentration, 
in accordance with the CEP for models of competition for a single resource \cite{tilman}.

The paper is organized as follows. In Section \ref{PR} we give some preliminary lemmas.
In Section \ref{wolkowicz} we show how the Lyapunov function of Wolkowicz and Lu \cite{wl} 
can be extended to (\ref{D=1}) and used to obtain Theorem \ref{ourthm1}. 
We show in this section that the result of \cite{sari1} for (\ref{eqsqi}),
which extends the result of \cite{wl} for (\ref{eqsxi}),
is a corollary of Theorem \ref{ourthm1}. We give also graphical interpretations of 
the conditions (\ref{h11}), (\ref{h31}) and (\ref{h21}) in Theorem \ref{ourthm1}.
In Section \ref{hsu}, we show how the Lyapunov function of Hsu \cite{hsu} can 
be extended to (\ref{D=1}) and used to obtain Theorem \ref{ourthm}, which is another  
global asymptotic stability result of $E_1^*$ for (\ref{D=1}).
Theorem \ref{ourthm} can be obtained also as a corollary of Theorem \ref{ourthm1} (see Proposition \ref{prop1}).
We show in this section that the result of \cite{hsu}  for 
(\ref{eqsxi}) with Monod functions (\ref{monod}), and the result of \cite{sari2} for (\ref{eqsqi})
are corollaries of Theorem \ref{ourthm}.  
In Section \ref{N=1} we discuss the single species case $N=1$.
In Section \ref{MLQ} we apply our result to the model with Monod growth functions (\ref{monod}) 
and linear yields.
In Section \ref{conclusions} we discuss some of the CEP results  
based on Lyapunov functions and we compare Theorem \ref{ourthm1} with the results of \cite{FH} 
based on a Bendixon-Dulac approach. 

\section{Preliminary results}\label{PR}
Let us prove first that we can we assume that $D=1$ and $S^0=1$ in (\ref{eqfpi}). 
Indeed, under the change of variables 
$$\overline{S}=\frac{S}{S^0},\qquad \overline{t}=Dt,\qquad 
\overline{p}_i(\overline{S})=\frac{p_i(S^0\overline{S})}{S^0D},\qquad 
\overline{f}_i(\overline{S})=\frac{f_i(S^0\overline{S})}{D},$$
equations (\ref{eqfpi}) take the form
$$
\begin{array}{lclcll} 
 \displaystyle \frac{d\overline{S}}{d\overline{t}}& =&
 \displaystyle \frac{1}{S^0D}\displaystyle \frac{d{S}}{d{t}}
&=&
1-\overline{S}-\displaystyle
\sum_{i=1}^N\overline{p}_i(\overline{S})x_i,&\\[3mm]
 \displaystyle 
\frac{d{x}_i}{d\overline{t}}
&=& 
 \displaystyle \frac{1}{D}\displaystyle \frac{d{x}_i}{d{t}}
&=&
\overline{f}_i(\overline{S})x_i,& i=1\cdots N.
\end{array}
$$
Dropping the bars, one obtains (\ref{D=1}).
Recall that $f_i(0)<0$, so that the concentration of the species $x_i$ 
is decreasing when the concentration of nutrient is too small. 
The smallest positive zero $S=\lambda_i$ of $f_i$ is the break-even concentration of 
the $i$th species $x_i$. We adopt the convention
$\lambda_i=\infty$ if $f_i(S)<0$ for all $S>0$. 
We need the following lemmas.
\begin{lemma}
The non-negative cone is invariant under the flow of (\ref{D=1}) 
and all solutions are defined and remain bounded for all $t\geq 0$.
\end{lemma}
This lemma is simply Theorem 4.1 in \cite{apw}.
\begin{lemma}\label{SppS1}
If for some species $x_i$, the inequality $(S-\lambda_i)f_i(S)>0$ is satisfied for all $0<S<1$, $S\neq\lambda_i$, 
then $S(t)< 1$ for all sufficiently large $t$ and all initial condition .
\end{lemma}
This lemma can be obtained using arguments similar to that given in the proofs of Lemma 2.9 in \cite{apw} and Lemma 2.1 in \cite{wl}.
\begin{lemma}\label{SppS0}
For all solutions of (\ref{D=1}), if $\lambda_i\geq 1$ then $x_i(t)\to 0$ as $t\to \infty$.
\end{lemma}
This lemma can be obtained using arguments similar to that given in the proofs of Lemma 4.2 in \cite{apw} and 
Lemma 2.2 in \cite{wl}.

\section{Extension of the Lyapunov function of Wolkowicz and Lu}\label{wolkowicz}
The Lyapunov function used by Wolkowicz and Lu \cite{wl} in the constant yields 
case (\ref{eqsxi}) is
\begin{equation}\label{LWL}
V_{WL}=
\frac{S^0-\lambda_1}{D_1}
\int_{\lambda_1}^S\frac{q_1(\sigma)-D_1}{S^0-\sigma}d\sigma+
\frac{1}{Y_1}\int_{x_{1}^*}^{x_1}\frac{\xi-x_{1}^*}{\xi}d\xi+
\sum_{i=2}^N\frac{c_i}{Y_i}x_i.
\end{equation}
with suitable constant $c_i>0$. 
Using the notations in (\ref{D=1}), and since $S^0$ was rescaled to 1, the numerator in the first integral of (\ref{LWL})
is simply equal to $\frac{f_1(\sigma)}{1-\sigma}$. Multiplying  (\ref{LWL}) by 
the constant 
$\frac{D_1}{1-\lambda_1}=\frac{Y_1}{x_1^*}$, gives the following
function 
\begin{equation}\label{lyapunov1}
V=
\int_{\lambda_1}^S\frac{f_1(\sigma)}{1-\sigma}d\sigma+
\frac{1}{x_{1}^*}\int_{x_{1}^*}^{x_1}\frac{\xi-x_{1}^*}{\xi}d\xi+
\sum_{i=2}^N\alpha_ix_i,
\end{equation}
where $\alpha_i$ are constants to be determined.
This is a Lyapunov function for (\ref{D=1}) which 
permits to prove Theorem \ref{ourthm1} as shown below.

{\em Proof}{\em{of Theorem \ref{ourthm1}}}.
From Lemmas \ref{SppS1} and \ref{SppS0} it follows that there is no loss of 
generality to assume that $\lambda_i<1$ for $i=1\cdots N$ and to restrict 
our attention to $0< S<1$.
Consider the function $V=V(S,x_1,\cdots,x_N)$ given by 
(\ref{lyapunov1})
where $\alpha_i$ are positive constants satisfying (\ref{h21}). 
The function $V$ is continuously 
differentiable for $0<S<1$ and $x_i>0$ and positive except at point $E_{1}^*$. 
The derivative of $V$ along the trajectories of (\ref{D=1}) is
$$
V'=\frac{f_1(S)}{1-S}S'+\frac{x_1-x_{1}^*}{x_{1}^*x_1}x_1'+
\sum_{i=2}^N\alpha_ix_i'.
$$
Since $x_{1}^*=P_1(\lambda_1)$ and using (\ref{D=1}), $V'$ is written 
$$
V'
=\frac{f_1(S)}{1-S}\left[1-S-\sum_{i=1}^Np_i(S)x_i\right]+
\frac{1}{P_1(\lambda_1)}[x_1-P_1(\lambda_1)]f_1(S)+
\sum_{i=2}^N\alpha_if_i(S)x_i.
$$
The terms $\frac{f_1(S)}{1-S}(1-S)$ and $-\frac{1}{P_1(\lambda_1)}P_1(\lambda_1)f_1(S)$ are canceled. Hence, using  (\ref{P1}),
$$
V'=x_1f_1(S)\left[\frac{1}{P_1(\lambda_1)}-\frac{1}{P_1(S)}\right]+
\sum_{i=2}^Nx_i\frac{\alpha_if_i(S)(1-S)-f_1(S)p_i(S)}{1-S}.
$$
Using (\ref{h11}) and (\ref{h31}), 
the first term of the above sum is non-positive
for $0<S<1$ and equals 0 if and only if $S=\lambda_1$ or $x_1=0$.  
Using (\ref{h21}), the second term 
is non-positive for $0<S<1$ and equals 0 if and only if $x_i=0$ for $i=2\cdots N$.
Hence $V'\leq 0$ and
$V'=0$ if and only if $x_i=0$ for $i=1\cdots N$ or 
$S=\lambda_1$ and $x_i=0$ for $i=2\cdots N$.
Using the Krasovskii-LaSalle extension theorem,
the $\omega$-limit set of the trajectory is $E_1^*$.
\medskip

Theorem \ref{ourthm1} was previously obtained in \cite{sari1},  
in the particular case when
the function $f_i$ 
has at most two positive zeros $\lambda_i$ and $\mu_i$, with 
$\lambda_i\leq \mu_i\leq +\infty$, such that
\begin{equation}\label{lambdamu}
f_i(S)<0\mbox{ if }S\notin[\lambda_i,\mu_i],\mbox{ and }
f_i(S)>0\mbox{ if }S\in]\lambda_i,\mu_i[,
\end{equation}
with the convention that $\mu_i= \infty$ if equation $f_i(S)=0$ 
has only one solution and $\lambda_i= \infty$ if 
it has no solution. This class of functions corresponds to the case when $f_i(S)=q_i(S)-D_i$ and
$$
q_i(S)<D_i\mbox{ if }S\notin[\lambda_i,\mu_i],\mbox{ and }
q_i(S)>D_i\mbox{ if }S\in]\lambda_i,\mu_i[.
$$
It was often considered in the literature \cite{bw,li,wl,wx}.
For this class of systems the main result in \cite{sari1} is
\begin{corollary}[Theorem 2.1 in \cite{sari1}]\label{ourthms}
Assume that 
\begin{equation}\label{hc11}
\lambda_1<\lambda_2\leq\cdots\leq\lambda_N,\mbox{ and }\lambda_1<1<\mu_1,
\end{equation}
\begin{equation}\label{hc31}
(S-\lambda_1)(P_1(S)-P_1(\lambda_1))<0,\mbox{ for }S\neq\lambda_1.
\end{equation}
There exist constants
$c_i>0$ for each $i\geq 2$ satisfying $\lambda_i<1$, such that
\begin{equation}\label{conditionci1}
\max_{0<S<\lambda_1}g_i(S)< c_i< \min_{\lambda_i<S<\rho_i}g_i(S),
\end{equation}
where 
$g_i(S)=\frac{1-\lambda_1}{p_1(\lambda_1)}\frac{f_1(S)p_i(S)}{f_i(S)(1-S)}$ and $\rho_i=\min(\mu_i,1)$.
Then the equilibrium $E_{1}^*$ is GAS for  (\ref{D=1}) with respect to the interior of the positive cone.
\end{corollary}
{\em Proof}
Assume that (\ref{hc11}), (\ref{hc31}) and (\ref{conditionci1}) hold. 
Let us prove that (\ref{h11}), (\ref{h31}) and (\ref{h21}) hold.
First, note that (\ref{hc31}) is the same as (\ref{h31}), 
and condition $\lambda_1<1<\mu_1$ in (\ref{hc11}) is equivalent 
to (\ref{h11}). If $\lambda_1<S<\lambda_i$ 
then $f_i(S)<0$ and $f_1(S)>0$ so that (\ref{h21})
is satisfied for any choice of $\alpha_i>0$.
Similarly if $\mu_i<1$ and 
$\mu_i<S<1$ then $f_i(S)<0$ 
and $f_1(S)>0$ so that (\ref{h21}) is satisfied for any choice of $\alpha_i>0$.
On the other hand, if 
$0<S<\lambda_1$ then $f_i(S)<0$ and, using $g_i(S)<c_i$ in (\ref{conditionci1}),
$$f_1(S)p_i(S)>c_i\frac{p_1(\lambda_1)}{1-\lambda_1}f_i(S)(1-S).$$
Finally, if $\lambda_i<S<\rho_i$, then $f_i(S)>0$ and, using $g_i(S)>c_i$ in (\ref{conditionci1}), 
$$ f_1(S)p_i(S)>c_i\frac{p_1(\lambda_1)}{1-\lambda_1}f_i(S)(1-S).$$ 
Thus (\ref{h21}) is satisfied for $\alpha_i=c_i\frac{p_1(\lambda_1)}{1-\lambda_1}$.
The result follows from Theorem \ref{ourthm1}.
\medskip 

Condition (\ref{h11}) means that $S=\lambda_1$ is the only zero of the growth function $f_1(S)$ for $0<S<1$.  
Condition (\ref{h31}) means that $S=\lambda_1$ is the only zero of the function $P_1(S)$ given by (\ref{P1}), 
for $0<S<1$.
The technical condition (\ref{h21}) is trivially satisfied in the single species $N=1$. 
Following \cite{sari2,wl} we give now a graphical interpretation of (\ref{h21}). For each $i\geq 2$ such that $\lambda_i<1$, consider 
the function
\begin{equation}\label{gi}
g_i(S)=\frac{f_i(S)}{f_1(S)}\frac{1-S}{p_i(S)}.
\end{equation}
The functions $g_i$ is defined on $(0,\lambda_1)\cup(\lambda_1,1]$.  It tends to $\pm\infty$ when $S$ tends $\lambda_1$. 
This function should not be confused with the function $g_i$ in Corollary \ref{ourthms} which is simply a multiple of 
the inverse of $g_i$. We use $g_i$ instead of its inverse, since the zeros of $f_i$ on $[0,1]$ are not known
as for the class of functions $f_i$ considered in Corollary~\ref{ourthms}.
Since $f_1(S)<0$ over $[0,\lambda_1)$ and $f_1(S)>0$ over $[\lambda_i,1]$, the condition (\ref{h21}) is equivalent to
\begin{equation}\label{h2graphic}
\min_{0<S<\lambda_1}g_i(S)> \frac{1}{\alpha_i}>
\max_{\lambda_i<S<1}g_i(S).
\end{equation}
Hence, the condition (\ref{h21}) in Theorem \ref{ourthm1} can be easily depicted  graphically:
plot simply the graph of $g_i(S)$ over $[0,1]$ and see if there is a gap between
 $\min_{0<S<\lambda_1}g_i(S)$ and $\max_{\lambda_i<S<1}g_i(S)$, see Fig. \ref{figF}.

It was shown in \cite{sari1} that the main result (Theorem 2.3) of \cite{wl} is a consequence of Corollary \ref{ourthms}. 
Hence, it is also a corollary of Theorem \ref{ourthm1}.

\section{Extension of the Lyapunov function of Hsu}\label{hsu}
The Lyapunov function $V_H$ used by Hsu \cite{hsu} 
in the Monod case of (\ref{eqsxi}), where the growth functions are of the form (\ref{monod}), 
is
\begin{equation}\label{HSU}
V_{H}=
\int_{\lambda_1}^S\frac{\sigma-\lambda_1}{\sigma}
d\sigma+c_1\int_{x_{1}^*}^{x_1}\frac{\xi-x_{1}^*}{\xi}d\xi+
\sum_{i=2}^Nc_ix_i,
\end{equation}
with constants 
$$c_i=\frac{1}{Y_i}\frac{a_i}{a_i-D_i},\quad i=\cdots N,\quad \mbox{and }\lambda_1=\frac{b_1D_1}{a_1-D_1}.$$
The function in the first integral of (\ref{HSU}) is simply given by
$\frac{S-\lambda_1}{S}
=c_1\frac{f_1(S)}{p_1(S)}$, where 
$$f_1(S)=\frac{a_1S}{b_1+S}-D_1,\qquad p_1(S)=\frac{1}{Y_1}\frac{a_1S}{b_1+S}.$$
Hence, multiplying (\ref{HSU}) by 
the constant $1/c_1$
gives the following function  
\begin{equation}\label{lyapunov}
V=
\int_{\lambda_1}^S\frac{f_1(\sigma)}{p_1(\sigma)}d\sigma+
\int_{x_{1}^*}^{x_1}\frac{\xi-x_{1}^*}{\xi}d\xi+
\sum_{i=2}^Nc_ix_i.
\end{equation}
where the constants $c_i/c_1$ in the last sum are simply denoted by $c_i$ to avoid unnecessary new notations.
Under some technical conditions, this function is a Lyapunov function for (\ref{D=1}) 
and permits to obtain the global asymptotic stability of 
the equilibrium point $E_{1}^*$ as stated in the following result.
\begin{theorem}\label{ourthm}
Assume that $\lambda_1<1$ and for all $0<S<1$,
\begin{equation}\label{h1}
(S-\lambda_1)f_1(S)>0,\mbox{ for }S\neq\lambda_1,
\end{equation}
\begin{equation}\label{h3}
(S-\lambda_1)(P_1(S)-P_1(\lambda_1))<0,\mbox{ for }S\neq\lambda_1.
\end{equation}
Assume that there exist constants
$c_i>0$ for each $i\geq 2$ satisfying $\lambda_i<1$, such that for all $0<S<1$,
\begin{equation}\label{h2}
f_1(S)p_i(S)>c_if_i(S)p_1(S).
\end{equation}
Then the equilibrium $E_{1}^*$ is GAS for (\ref{D=1}) with respect to the interior of the positive cone.
\end{theorem}
{\em Proof}
From Lemmas \ref{SppS1} and \ref{SppS0} it follows that there is no loss of generality to assume that $\lambda_i<1$ for $i=1\cdots N$ and to restrict 
our attention to $0< S<1$.
Consider the function $V=V(S,x_1,\cdots,x_N)$ given by 
(\ref{lyapunov})
where $c_i$ are positive constants satisfying (\ref{h2}). 
The function $V$ is continuously 
differentiable in the positive cone and positive except at point $E_{1}^*$, where it is equal to 0. 
The derivative of $V$ along the trajectories of (\ref{D=1}) is
$$
V'=\frac{f_1(S)}{p_1(S)}S'+\frac{x_1-x_{1}^*}{x_1}x_1'+
\sum_{i=2}^Nc_ix_i'.
$$
Since $x_{1}^*=P_1(\lambda_1)$ and, using (\ref{D=1}), $V'$ is written
$$
V'=\frac{f_1(S)}{p_1(S)}\left[1-S-\sum_{i=1}^Np_i(S)x_i\right]+
[x_1-P_1(\lambda_1)]f_1(S)+
\sum_{i=2}^Nc_if_i(S)x_i.
$$
The terms $-\frac{f_1(S)}{p_1(S)}p_1(S)x_1$ and $x_1f_1(S)$ are canceled. Therefore, using (\ref{P1}),
$$
V'=\displaystyle f_1(S)\left[P_1(S)-P_1(\lambda_1)\right]+
\sum_{i=2}^Nx_i\frac{c_if_i(S)p_1(S)-f_1(S)p_i(S)}{p_1(S)}.
$$
Using (\ref{h1}) and (\ref{h3}), 
the first term of the above sum is non-positive
for $0<S<1$ and equals $0$ if and only if $S=\lambda_1$. 
Using (\ref{h2}), the second term 
is non-positive
for $0<S<1$ and equals $0$ if and only if $x_i=0$ for $i=2\cdots N$.
Hence $V'\leq 0$ and
$V'=0$ if and only if $S=\lambda_1$ and $x_i=0$ for $i=2\cdots N$.
By the Krasovskii-LaSalle extension theorem, the $\omega$-limit set of the trajectory is $E_1^*$.
\medskip

Notice that, as in Lemma \ref{lemmaCEP}, the condition (\ref{h2}) implies that $\lambda_1<\lambda_i$ for all $i\geq 2$. 
Indeed, if $f_i(S)>0$ for some $S\leq\lambda_1$, then $f_1(S)\leq 0$, so that the inequality (\ref{h2}) 
is violated. Thus the winning species $x_1$ in Theorem \ref{ourthm} has the lowest break-even concentration.
Actually Theorem \ref{ourthm} is a consequence of Theorem \ref{ourthm1}.

\begin{proposition}\label{prop1}
Theorem \ref{ourthm} is a corollary of Theorem \ref{ourthm1}.
\end{proposition}
{\em Proof}
Assume that (\ref{h1}), (\ref{h3}) and (\ref{h2}) hold. Notice that (\ref{h1}) is the same as (\ref{h11})
and (\ref{h3}) is the same as (\ref{h31}). Let us prove that (\ref{h21}) holds.
If $f_i(S)<0$ and $f_1(S)>0$ (which occurs if $\lambda_1<S<\lambda_i$ and may occur also for $\lambda_i<S<1$), 
then (\ref{h21}) holds for any choice of $\alpha_i>0$.
If $0<S<\lambda_1$ then, by (\ref{h3}), $P_1(S)>P_1(\lambda_1)$ and, since $f_i(S)<0$, 
$\frac{f_i(S)}{P_1(S)}>\frac{f_i(S)}{P_1(\lambda_1)}$.
Finally, if $\lambda_i<S<1$ and $f_i(S)>0$ then, by (\ref{h3}), $P_1(S)<P_1(\lambda_1)$, and hence,
$\frac{f_i(S)}{P_1(S)}>\frac{f_i(S)}{P_1(\lambda_1)}$. 
Therefore, in both cases $\lambda_i<S<1$ and $0<S<\lambda_1$,
$$\frac{f_i(S)}{P_1(S)}>\frac{f_i(S)}{P_1(\lambda_1)}.$$ 
Thus, using (\ref{h2}), 
$$f_1(S)p_i(S)>c_i{f_i(S)p_1(S)}=c_i\frac{f_i(S)}{P_1(S)}(1-S)>c_i\frac{f_i(S)}{P_1(\lambda_1)}(1-S).$$
Thus, (\ref{h21}) holds for $\alpha_i=\frac{c_i}{P_1(\lambda_1)}$.
Hence, (\ref{h3}) and (\ref{h2}) imply (\ref{h21}). This ends the proof. 
\medskip

Theorem \ref{ourthm} recovers the classical case Monod case \cite{hsu}. Indeed,
consider the particular case of (\ref{eqsxi}), when the growth functions $q_i(S)$ are given by (\ref{monod}).
System (\ref{eqsxi}), with $D=1$ and $S^0=1$, takes the form
\begin{equation}\label{eqsxmm}
\begin{array}{l} 
 \displaystyle S' =1-S-\sum_{i=1}^N\frac{a_iS}{b_i+S}\frac{x_i}{Y_i},\\[3mm]
 \displaystyle x'_i= \left[\frac{a_iS}{b_i+S} - D_i\right]x_i,\qquad\qquad i=1\cdots N.
\end{array}
\end{equation}
We consider the case where, for all $i=1\cdots N$, $a_i>D_i$.
The break-even concentrations are
\begin{equation}\label{bec}
\lambda_i=\frac{b_iD_i}{a_i-D_i}.
\end{equation}
\begin{corollary}[Theorem 3.3 in \cite{hsu}]\label{thmhsu}
Assume that 
\begin{equation}\label{hcm1}
\lambda_1<\lambda_2\leq\cdots\leq\lambda_N,\qquad \lambda_1<1.
\end{equation}
Then the equilibrium $E_{1}^*$ is GAS
for (\ref{eqsxmm}) with respect to the interior of the positive cone.
\end{corollary}
{\em Proof}
Assume that (\ref{hcm1}) holds. Let us prove that (\ref{h1}), (\ref{h1}) and (\ref{h1}) hold.
Since $f_1(S)=q_1(S)-D_1$ is increasing, the function $f_1(S)$ changes sign only at $S=\lambda_1$ and hence,
(\ref{h1}) is satisfied. Since
$$P_1(S)=Y_1(1-S)\frac{b_1+S}{a_1S}
\qquad
\mbox{ and }
\qquad
P'_1(S)=-Y_1\frac{S^2+b_1}{a_1S^2}<0,
$$
the function $P_1(S)$ changes sign only at $S=\lambda_1$ and hence (\ref{h3}) is satisfied.
Condition  (\ref{h2}) is
$$
\frac{(a_1-D_1)S-b_1D_1}{b_1+S}\frac{1}{Y_i}\frac{a_iS}{b_i+S}>
c_i\frac{(a_i-D_i)S-b_iD_i}{b_i+S}\frac{1}{Y_1}\frac{a_1S}{b_1+S},\qquad i\geq 2.
$$
After simplification by $\frac{S}{(b_1+S)(b_i+S)}$, this condition is equivalent to 
\begin{equation}\label{h2m}
(a_1-D_1)\frac{a_i}{Y_i}(S-\lambda_1)
>c_i(a_i-D_i)\frac{a_1}{Y_1}(S-\lambda_i)\qquad i\geq 2,
\end{equation}
which is satisfied for $c_i=\frac{(a_1-D_1)a_iY_1}{(a_i-D_i)a_1Y_i}$. Indeed, 
for this choice of the constants $c_i$, (\ref{h2m}) is simply
$$S-\lambda_1>S-\lambda_i\Longleftrightarrow \lambda_1<\lambda_i,\qquad i\geq 2,$$
which is the same as (\ref{hcm1}).
Thus  (\ref{h2}) is satisfied.
The global asymptotic stability of 
$E_{1}^*$ follows from Theorem \ref{ourthm}.
\medskip

Theorem \ref{ourthm}  was previously obtained in \cite{sari2}, 
in the particular case when the function $f_i$ satisfies (\ref{lambdamu}).
For this class of systems the main result in \cite{sari2} is
\begin{corollary}[Theorem 2 in \cite{sari2}]\label{ourthmsm}
Assume that 
\begin{equation}\label{hc1}
\lambda_1<\lambda_2\leq\cdots\leq\lambda_N,\mbox{ and }\lambda_1<1<\mu_1,
\end{equation}
\begin{equation}\label{hc3}
(S-\lambda_1)(P_1(S)-P_1(\lambda_1))<0,\mbox{ for }S\neq\lambda_1.
\end{equation}
There exist constants
$\alpha_i>0$ for each $i\geq 2$ satisfying $\lambda_i<1$, such that
\begin{equation}\label{conditionci}
\max_{0<S<\lambda_1}g_i(S)< c_i< \min_{\lambda_i<S<\rho_i}g_i(S),
\end{equation}
where 
$g_i(S)=\frac{f_1(S)p_i(S)}{f_i(S)p_1(S)}$ and $\rho_i=\min(\mu_i,1)$.
Then the equilibrium $E_{1}^*$ is GAS
for (\ref{eqfpi}) with respect to the interior of the positive cone.
\end{corollary}
{\em Proof}
First, note that (\ref{hc3}) is the same as (\ref{h3}), 
and condition $\lambda_1<1<\mu_1$ in (\ref{hc1}) is equivalent 
to (\ref{h1}). If $S<\lambda_i$ 
then $f_i(S)<0$ and $f_1(S)>0$ so that (\ref{h2}) 
is satisfied for any choice of $c_i>0$.
Similarly if $\mu_i<1$ and 
$\mu_i<S<1$ then $f_i(S)<0$ 
and $f_1(S)>0$ so that (\ref{h2}) 
is satisfied for any choice of $c_i>0$.
On the other hand, if 
$0<S<\lambda_1$ then $f_i(S)<0$ and, using $g_i(S)<c_i$ in (\ref{conditionci}), 
$f_1(S)p_i(S)<c_if_i(S)p_1(S)$. 
Finally, if $\lambda_i<S<\rho_i$, then $f_i(S)>0$ and, using $g_i(S)>c_i$ in (\ref{conditionci}), 
$f_1(S)p_i(S)<c_if_i(S)p_1(S)$.
Thus (\ref{h2}) is satisfied.
The result follows from Theorem \ref{ourthm}.

\section{Single species}\label{N=1}
In the case $N=1$, (\ref{D=1}) takes the form
\begin{equation}\label{eqsx}
\begin{array}{l}
S' = 1-S-xp(S),\\
x'= f(S)x.
\end{array}
\end{equation}
Let $S=\lambda$ be the smallest positive value of $S$ such that $f(S)=0$ and $x^*=P(\lambda)$ with $P(S)$ 
defined by $P(S)=\frac{1-S}{p(S)}$ as in (\ref{P1}). 
If $\lambda<1$, then  $E^*=(\lambda,x^*)$ is a positive equilibrium. 
Assume that $f'(\lambda)>0$ and $P'(\lambda)<0$, so that $E^*$ is locally asymptotically stable. 
We consider the global asymptotic stability of $E^*$.

\begin{corollary}[Theorem 2.11 in \cite{apw} or Lemma 2.3 in \cite{PW}]\label{thmN=1}
Assume that $\lambda<1$ and for all $0<S<1$,
\begin{equation}\label{h1N=1}
(S-\lambda)f(S)>0,\mbox{ for }S\neq\lambda,
\end{equation}
\begin{equation}\label{h3N=1}
(S-\lambda)(P(S)-P(\lambda))<0,\mbox{ for }S\neq\lambda.
\end{equation}
Then the equilibrium $E^*$ is GAS
for (\ref{eqsx}) with respect to the interior of the positive cone.
\end{corollary}
{\em Proof}
Notice that (\ref{h1N=1}) is the same as (\ref{h11}) or (\ref{h1}) 
and (\ref{h3N=1}) is the same as (\ref{h31}) or (\ref{h3}). 
Since for $N=1$,  condition (\ref{h21}) in Theorem \ref{ourthm1} 
or condition (\ref{h2}) in Theorem \ref{ourthm}
is trivially satisfied, the result is a corollary of 
Theorem \ref{ourthm1} or Theorem \ref{ourthm}.
\medskip

Corollary \ref{thmN=1} was obtained  by Arino, 
Pilyugin and Wolkowicz (see \cite{apw}, Theorem 2.11). 
Using the Lyapunov function
\begin{equation}\label{APW}
V_{APW}=\frac{1-\lambda}{p(\lambda)}
\int_{\lambda}^S\frac{f(\sigma)}{1-\sigma}d\sigma+
\int_{x^*}^{x}\frac{\xi-x^*}{\xi}d\xi,
\end{equation}
these authors proved that if 
\begin{equation}\label{condAPW}
1-\frac{p(S)(1-\lambda)}{p(\lambda)(1-S)}\mbox{ has exactly one sign change for }0<S<1
\end{equation}
then $E^*$ is globally asymptotically stable. 
Notice that (\ref{condAPW}) is equivalent to (\ref{h3N=1}). 
In the single species case, our Lyapunov function (\ref{lyapunov1}), 
used in the proof of Theorem \ref{ourthm1}, 
reduces (up to a constant) to the Lyapunov function $V_{APW}$ considered in \cite{apw}.
Corollary \ref{thmN=1} was obtained also by  
Pilyugin and Waltman (see \cite{PW}, Lemma 2.3). Using the Lyapunov function
\begin{equation}\label{LPW}
V_{PW}=
\int_{\lambda}^S\frac{f(\sigma)}{p(\sigma)}d\sigma+
\int_{x^*}^{x}\frac{\xi-x^*}{\xi}d\xi,
\end{equation}
these authors proved that if 
\begin{equation}\label{condPW}
S=\lambda\mbox{ is the only zero of }R(S)=1-S-x^*p(S)
\end{equation} 
then $E^*$ is globally asymptotically stable. 
Notice that (\ref{condPW}) is equivalent to (\ref{h3N=1}). 
In the single species case, our Lyapunov function (\ref{lyapunov}), 
used in the proof of Theorem \ref{ourthm}, 
reduces to the Lyapunov function $V_{PW}$ considered in \cite{PW}.

\begin{figure}[ht]
\setlength{\unitlength}{1.0cm}
\begin{center}
\begin{picture}(13,3.5)(0,0)
\put(-0.5,0){\rotatebox{0}{\includegraphics[scale=0.16]{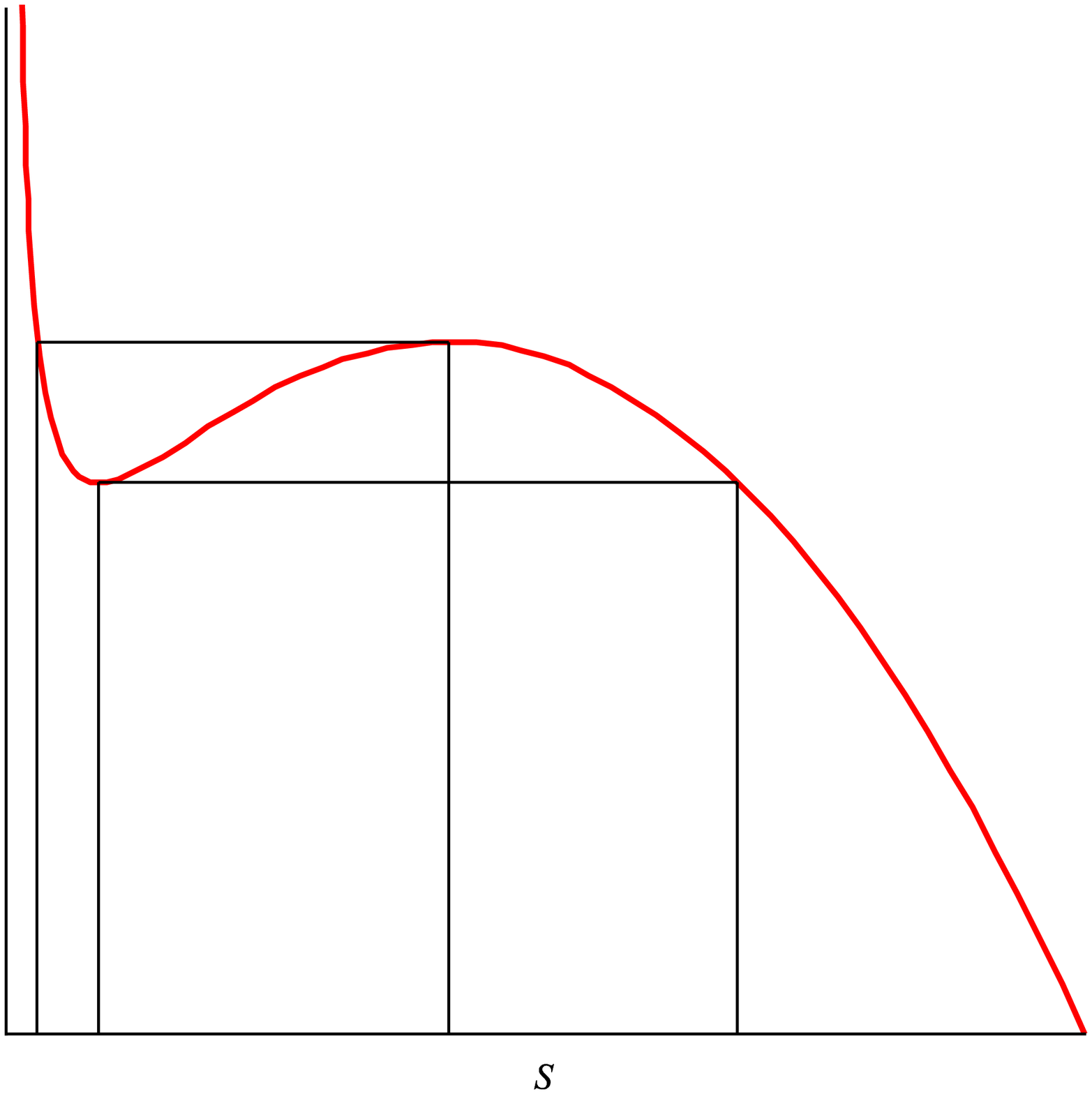}}}
\put(3,0){\rotatebox{0}{\includegraphics[scale=0.16]{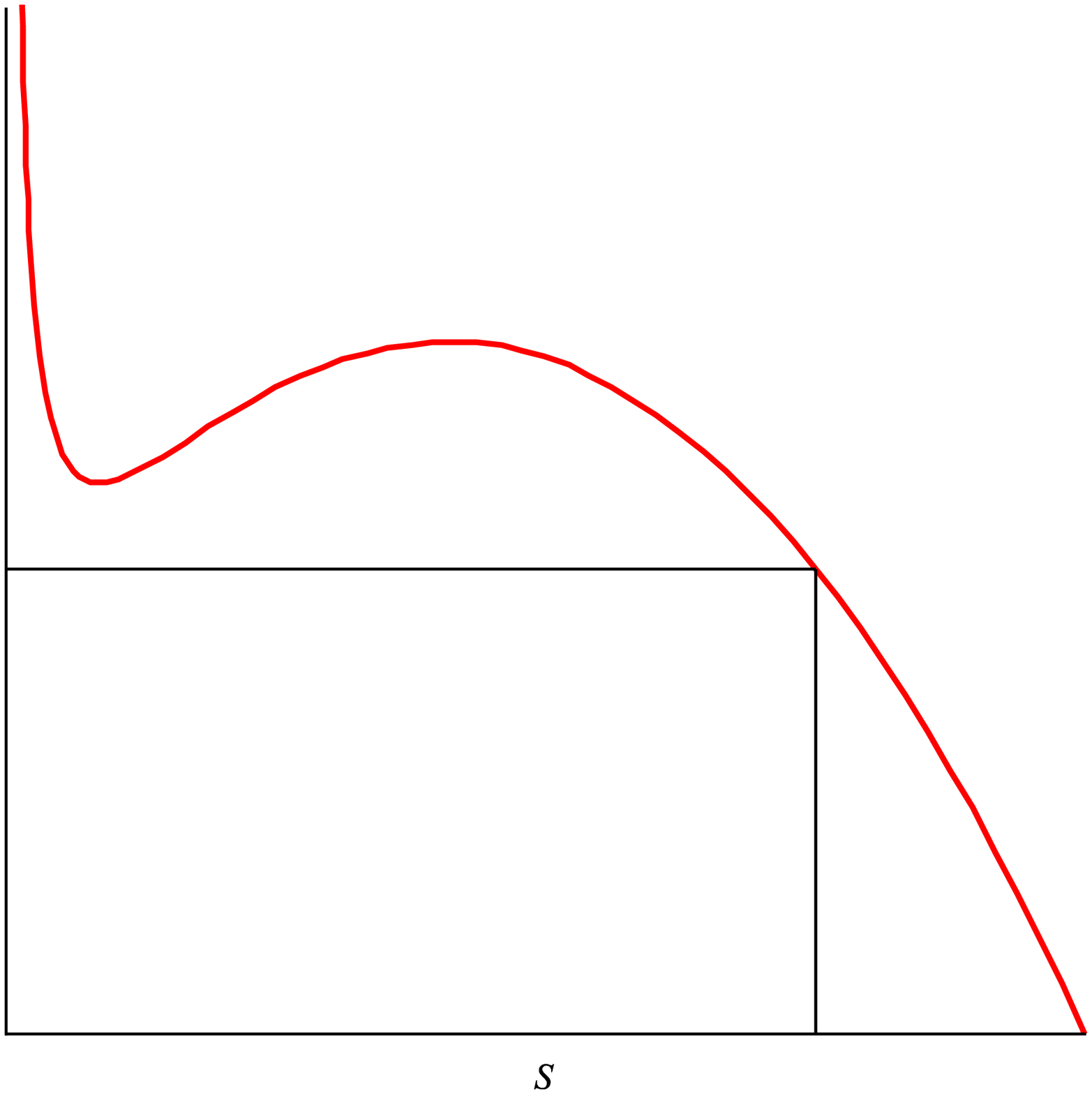}}}
\put(6.5,0){\rotatebox{0}{\includegraphics[scale=0.16]{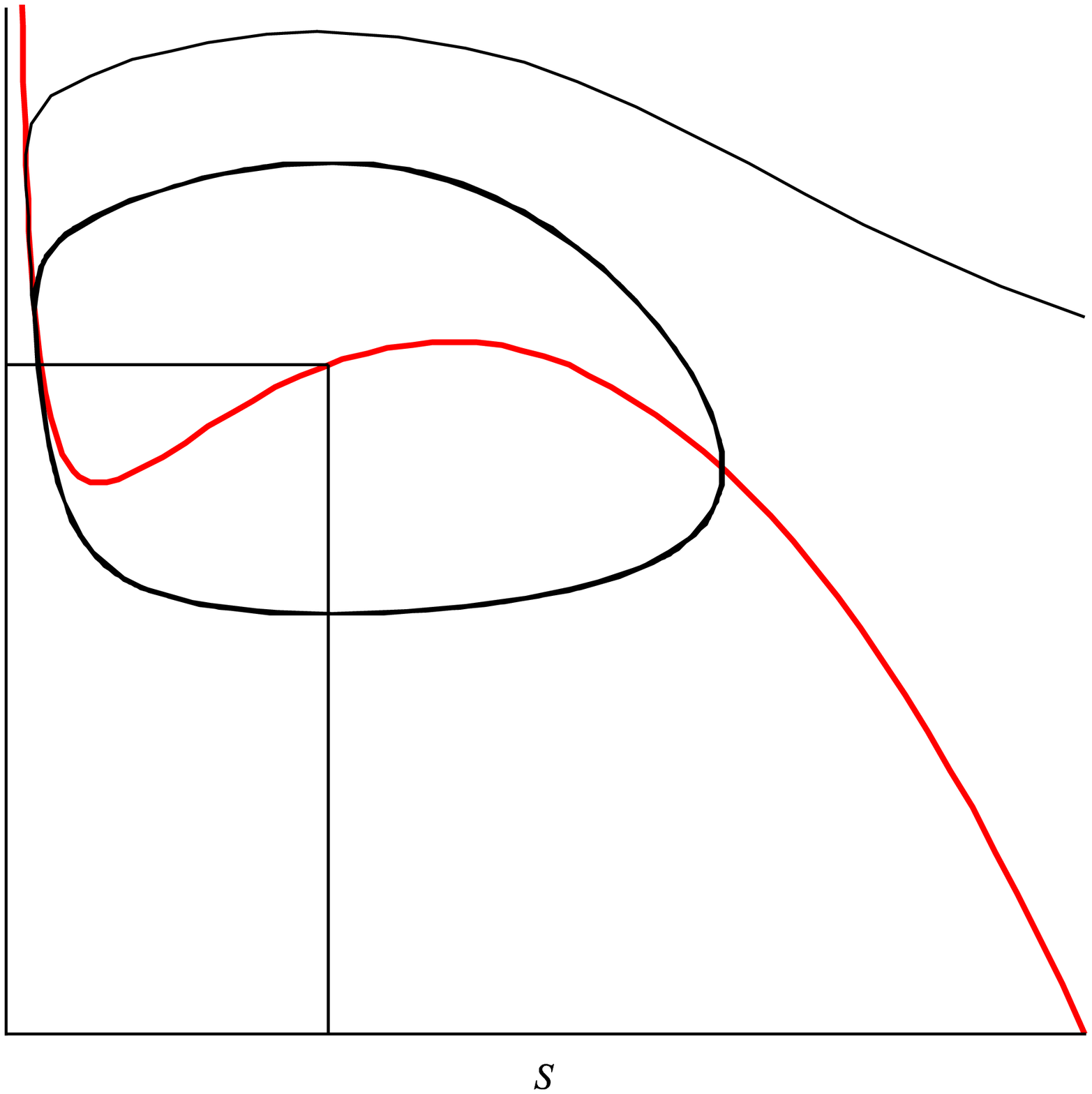}}}
\put(10,0){\rotatebox{0}{\includegraphics[scale=0.16]{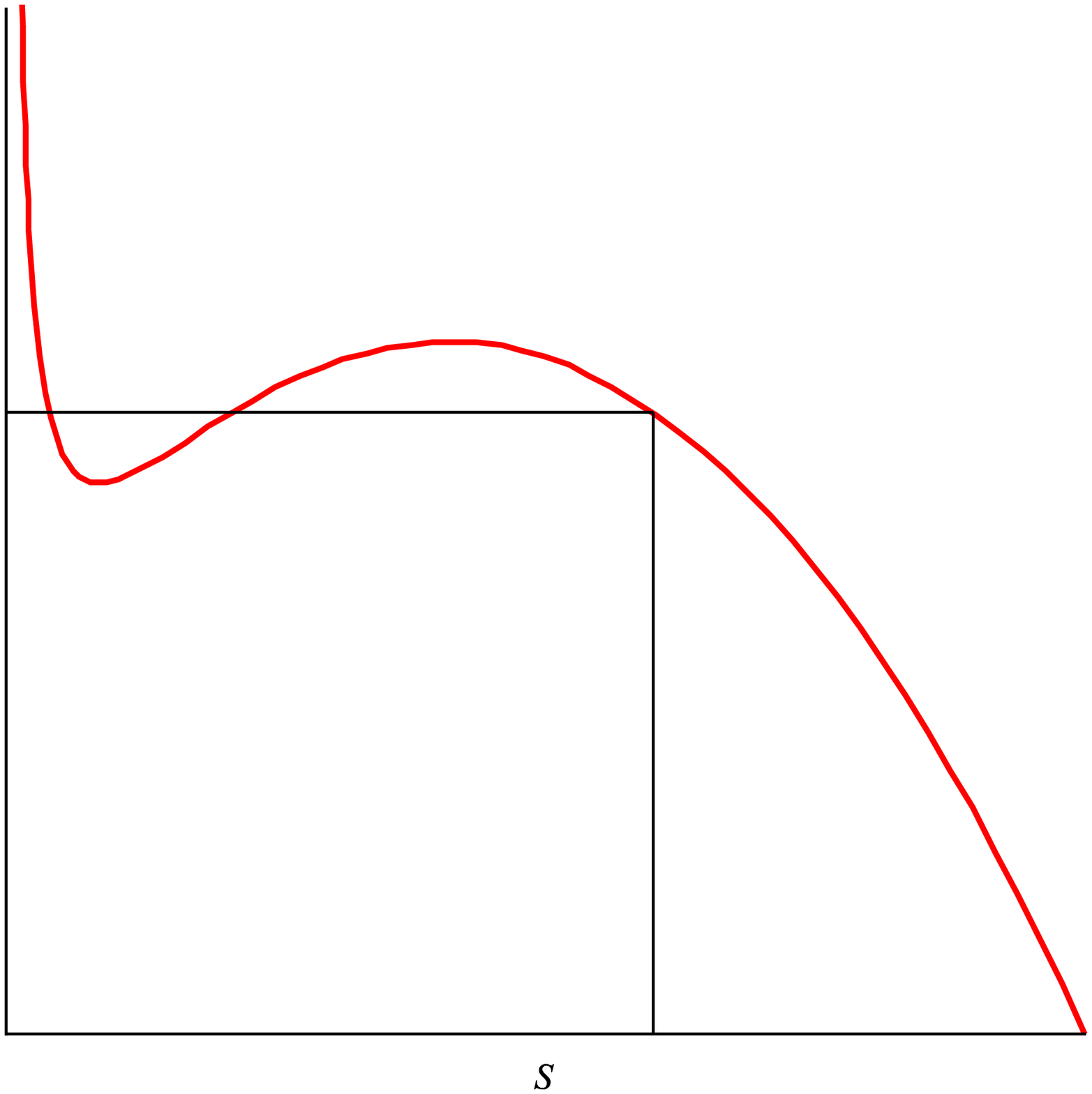}}}
\put(-0.5,0){{\small $S_1$}}
\put(-0.2,0){{\small $S_2$}}
\put(0.7,0){{\small $S_3$}}
\put(1.5,0){{\small $S_4$}}
\put(-0.3,2.9){{\small $x=P(S)$}}
\put(3.5,-0.3){{\small $\lambda<S_1$ or $\lambda>S_4$}}
\put(7.2,-0.3){{\small $S_2<\lambda<S_3$}}
\put(10,-0.3){{\small $S_1<\lambda<S_2$ or $S_3<\lambda<S_4$}}
\put(5.2,-0.05){{\small $\lambda$}}
\put(7.4,-0.05){{\small $\lambda$}}
\put(11.8,-0.08){{\small $\lambda$}}
\put(2.8,1.4){{\small $x^*$}}
\put(6.3,2){{\small $x^*$}}
\put(9.8,1.9){{\small $x^*$}}
\put(5.4,1.5){{\small $E^*$}}
\put(7.3,2.15){{\small $E^*$}}
\put(11.9,2){{\small $E^*$}}
\end{picture}
\end{center}
\caption{The graph of the function $x=P(S)$ showing the values $S_1$, $S_2$, $S_3$ and $S_4$.} \label{fig1}
\end{figure}

\begin{figure}[ht]
\setlength{\unitlength}{1.0cm}
\begin{center}
\begin{picture}(13,3.8)(0,0)
\put(0,0){\rotatebox{0}{\includegraphics[scale=0.2]{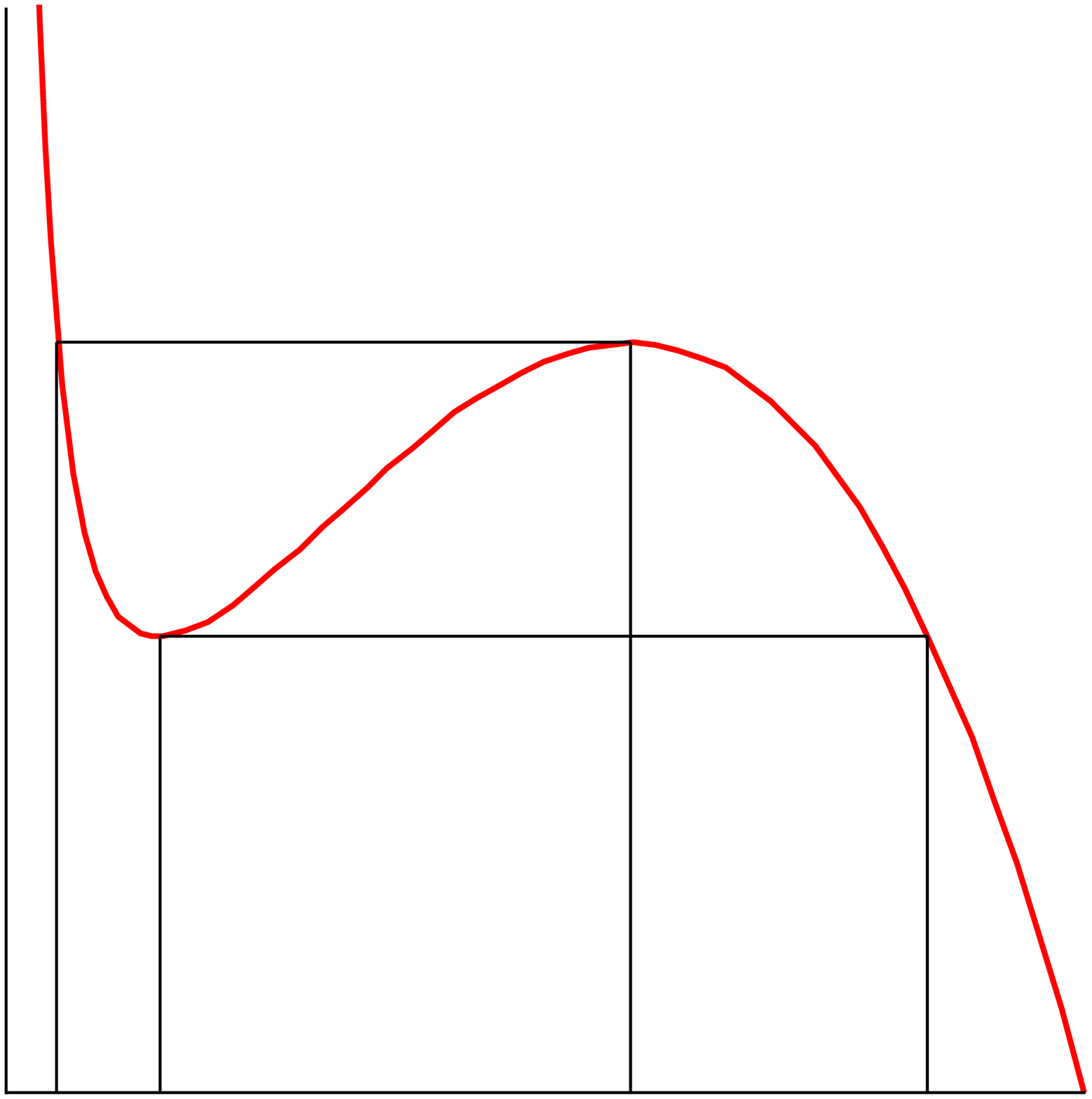}}}
\put(4.5,0){\rotatebox{0}{\includegraphics[scale=0.2]{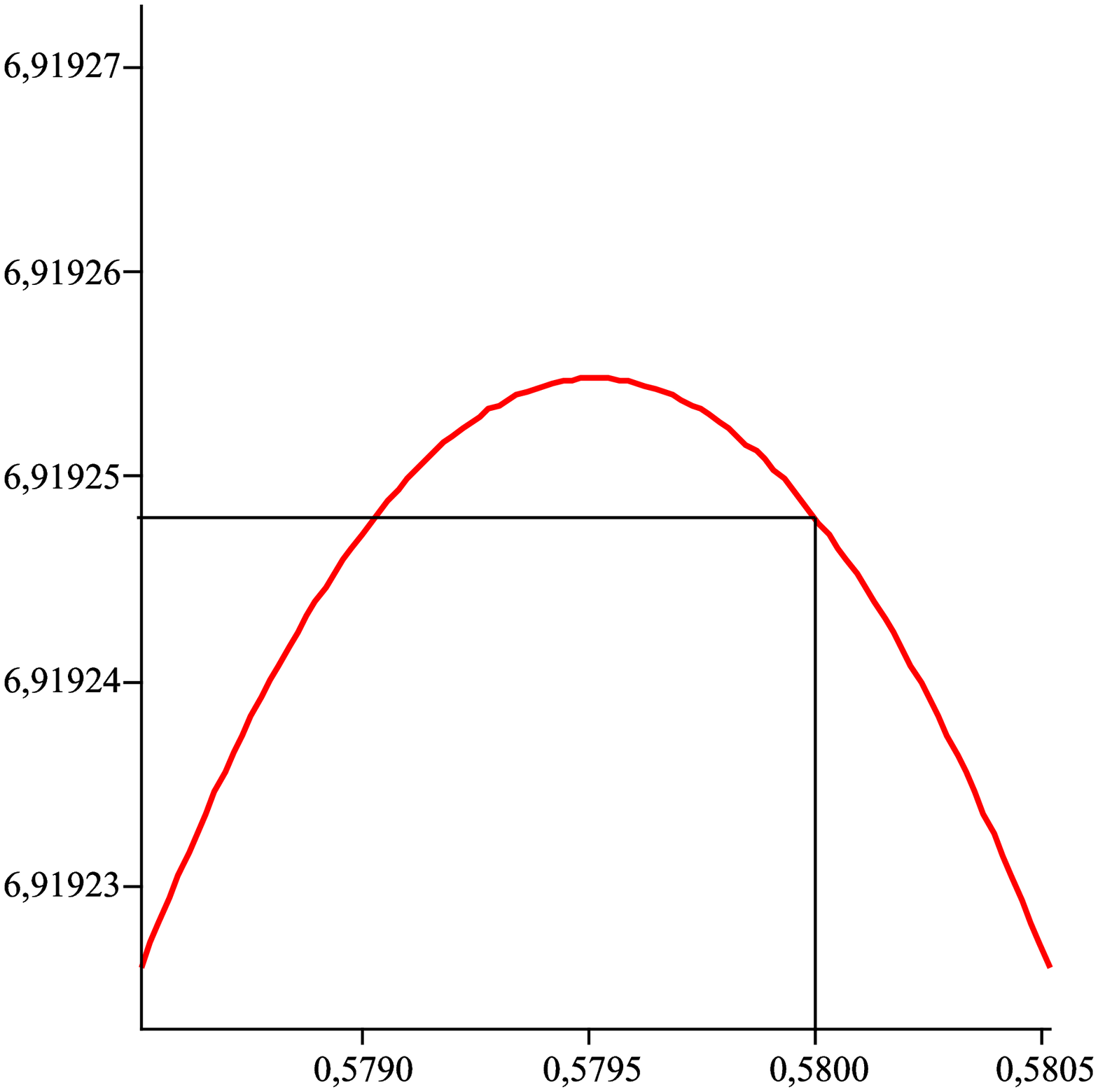}}}
\put(9,0){\rotatebox{0}{\includegraphics[scale=0.2]{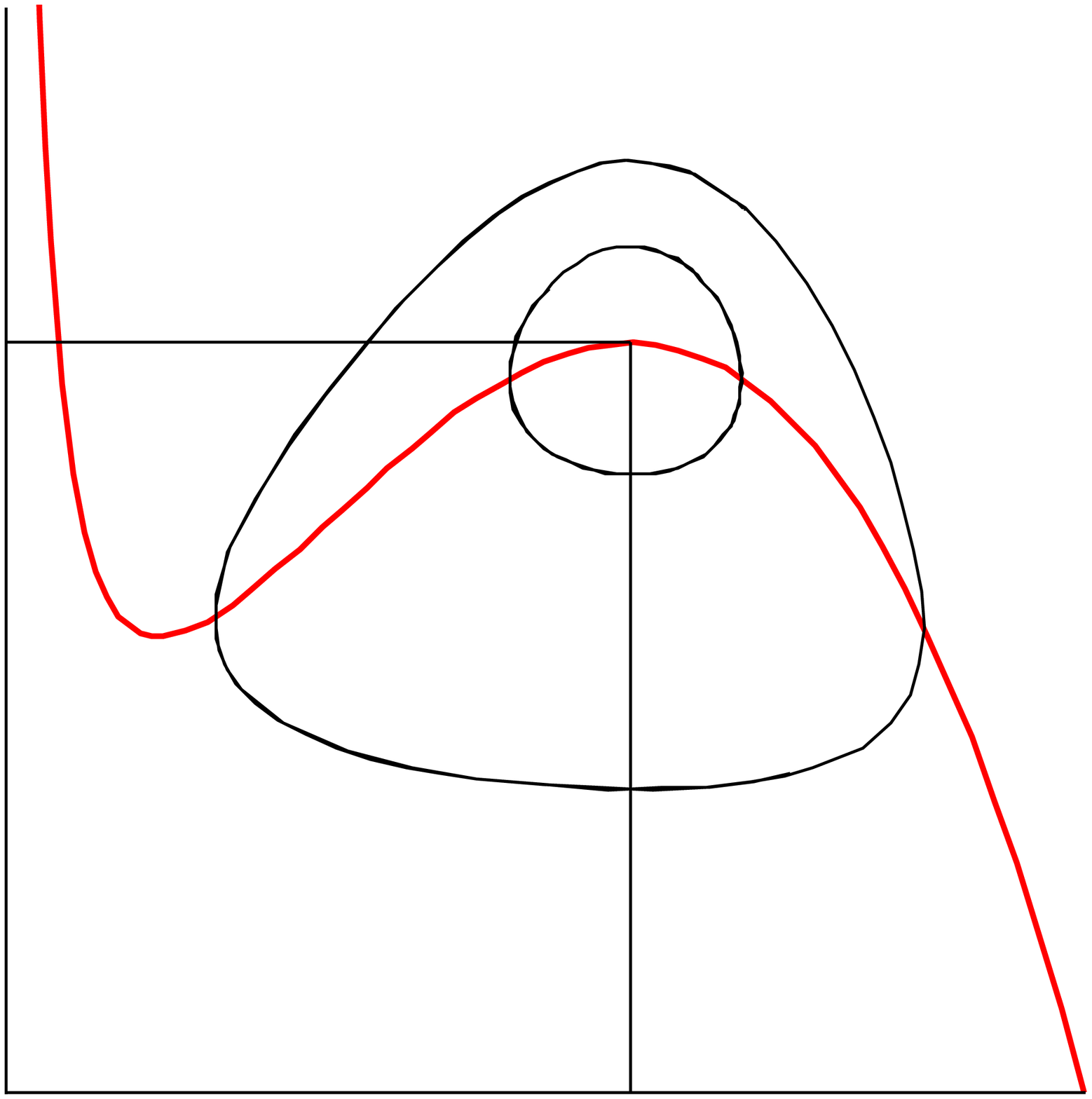}}}
\put(0.2,-0.2){{\small $S_1$}}
\put(0.55,-0.2){{\small $S_2$}}
\put(2.2,-0.2){{\small $S_3$}}
\put(3.2,-0.2){{\small $S_4$}}
\put(0.3,3.2){{\small $x=P(S)$}}
\put(6,2.6){{\small $x=P(S)$}}
\put(7.4,2.1){{\small $E^*$}}
\end{picture}
\end{center}
\caption{
If $\lambda<S_1$ or $\lambda>S_4$, using Corollary \ref{thmN=1}, the equilibrium point $E^*$ is GAS. 
If $S_1<\lambda<S_2$ or $S_3<\lambda<S_4$
the condition (\ref{h3N=1}) does not hold and Corollary \ref{thmN=1} cannot be applied. 
Actually if $D=1$ the system has two limit cycles. On the center of the figure an enlargement 
of the graph shows that the equilibrium point $E^*=(\lambda,x^*)$
lies on a decreasing branch of the graph of the function $x=P(S)$. 
For the parameters values $a=2$, $b=0.58$ and $c=46$: 
$S_1\simeq0.048$, $S_2\simeq0.143$, $S_3\simeq0.579$, $S_4\simeq0.855$ 
and $\lambda=0.58$. Hence $S_3<\lambda<S_4$. 
} \label{fig2}
\end{figure}

Notice that 
the isoclines $S'=0$ and $x'=0$ of (\ref{eqsx}) are given by
$$S'=0\Longleftrightarrow x=P(S),$$
$$x'=0\Longleftrightarrow x=0 \mbox{ or } S=\lambda.$$
If the vertical line $S=\lambda$ intersects the curve $x=P(S)$ on an increasing arc, 
then, from Lemma \ref{LAS}, the intersection is an unstable equilibrium point $E^*$. 
Using Poincar\'e-Bendixon theory we can show that the system has at least a periodic 
orbit surrounding the equilibrium.
Otherwise, if the vertical line $S=\lambda$ intersects the curve $x=P(S)$ on an decreasing arc, 
then, from Lemma \ref{LAS}, $E^*$
is locally asymptotically stable. 
The condition (\ref{h3N=1}) has the following graphical interpretation: 
if $E^*$ is the only intersection of the isocline $x=P(S)$ 
with the horizontal line $x=x^*$ then $E^*$ is GAS. 
For instance, in the situation depicted on Fig. \ref{fig2}, the function $x=P(S)$ has two critical points
$S=S_2$ and $S=S_3$. Let $S_1$ and $S_2$ defined by $P(S_1)=P(S_3)$ and $P(S_4)=P(S_2)$ respectively. Then:

Case 1. If $0<\lambda<S_1$ or $S_4<\lambda<1$ then $E^*$ is the only intersection of the isocline $x=P(S)$ 
with the horizontal line $x=x^*$. Thus, using Corollary \ref{thmN=1}, the equilibrium $E^*$ is GAS.

Case 2. If $S_2<\lambda<S_3$ then, using Lemma \ref{LAS}, the equilibrium $E^*$ is unstable. The 
system admits at least one limit cycle.

Case 3. If $S_1<\lambda<S_2$ or $S_3<\lambda<S_4$ then, using Lemma \ref{LAS}, the equilibrium $E^*$ 
is locally asymptotically stable. 
The horizontal line  $x=x^*$ has three intersections with $x=P(S)$.
Since (\ref{h3N=1}) does not hold, we cannot conclude if the equilibrium is GAS or not.
\medskip

We illustrate the third case by an example taken from  \cite{PW}. Consider  (\ref{eqsx}) with 
$$p(S)=\frac{q(S)}{y(S)}\mbox{ and }f(S)=q(S)-D,\mbox{ where } q(S)=\frac{aS}{b+S}\mbox{ and } y(S)=1+cS^2,$$
corresponding to Monod growth function and quadratic yield. 
For the parameter values given in the caption of Fig. \ref{fig2}, and  $D>D_3$ and close to $D_3$, 
the equilibrium point $E^*$ is exponentially stable and it is surrounded by two limit cycles. Actually, the limit cycle which exists for all $D_2<D<D_3$ disappears
for some critical $D_c>D_3$ through a subcritical Hopf bifurcation. 
For more details and explanations the reader is referred to \cite{PW}.

The case $N=1$ of a single species can also be investigated with the Bendixon-Dulac criterion. 
As shown in \cite{FH}, using the variable $y=\log(x)$, (\ref{eqsx}) is written 
\begin{equation}\label{eqsx1}
\begin{array}{l}
S' = 1-S-e^yp(S),\\
y'= f(S),
\end{array}
\end{equation}
with resulting divergence
$${\rm div}=-1-e^yp'(S).$$
If $p'(S)<0$ then the divergence is negative and no periodic solution can exist. 
Poincar\'e-Bendixon theorem shows that convergence to the equilibrium $E^*$ ensues.
Thus $E^*$ is GAS under the conditions $p'>0$ and (\ref{h1N=1}). 
This result is a consequence of Corollary \ref{thmN=1},  
since the condition $p'(S)$ implies $P'(S)<0$ for $0<S<1$ 
and hence, (\ref{h3N=1}) holds. 
However, the condition (\ref{h3N=1}) in Corollary \ref{thmN=1} can accept slightly negative $p'$, since
$P'(S)<0$, for $0<S<1$, is equivalent to $p'(S)>\frac{-p(S)}{1-S}$, for $0<S<1$.

\section{Monod growth functions and linear yields}\label{MLQ}
Models with linear yields were biologically motivated by \cite{ALLR,CT,CWT} 
who noticed the existence of limit cycles for some values of the parameters. 
The rigorous mathematical study was given in \cite{PW}. 
Consider the particular case of (\ref{eqsqi}), where the growth functions $q_i(S)$
are given by (\ref{monod}), and the yields $y_i(S)=p_i(S)/f_i(S)$ are linear
$$
y_i(S)=Y_i(1+c_iS)
$$
where $Y_i>0$ and $c_i\geq 0$. 
System (\ref{eqsqi}), with $D=1$ and $S^0=1$, takes the form
\begin{equation}\label{eqMonodLQ}
\begin{array}{l}
 \displaystyle S' =1-S-\sum_{i=1}^N\frac{a_iS}{b_i+S}\frac{x_i}{Y_i(1+c_iS)},\\
 \displaystyle x'_i= \left[\frac{a_iS}{b_i+S} - D_i\right]x_i,\qquad\qquad i=1\cdots N.
\end{array}
\end{equation}
The break-even concentrations $\lambda_i$ are given by (\ref{bec}).
In this section we give analytical conditions on the parameters of (\ref{eqMonodLQ})
so that conditions (\ref{h11}), (\ref{h31}) and (\ref{h21}) are satisfied and 
Theorem \ref{ourthm1} can be applied. 
\begin{figure}[ht]
\setlength{\unitlength}{1.0cm}
\begin{center}
\begin{picture}(13,3.5)(0,0)
\put(-0.5,0){\rotatebox{0}{\includegraphics[scale=0.16]{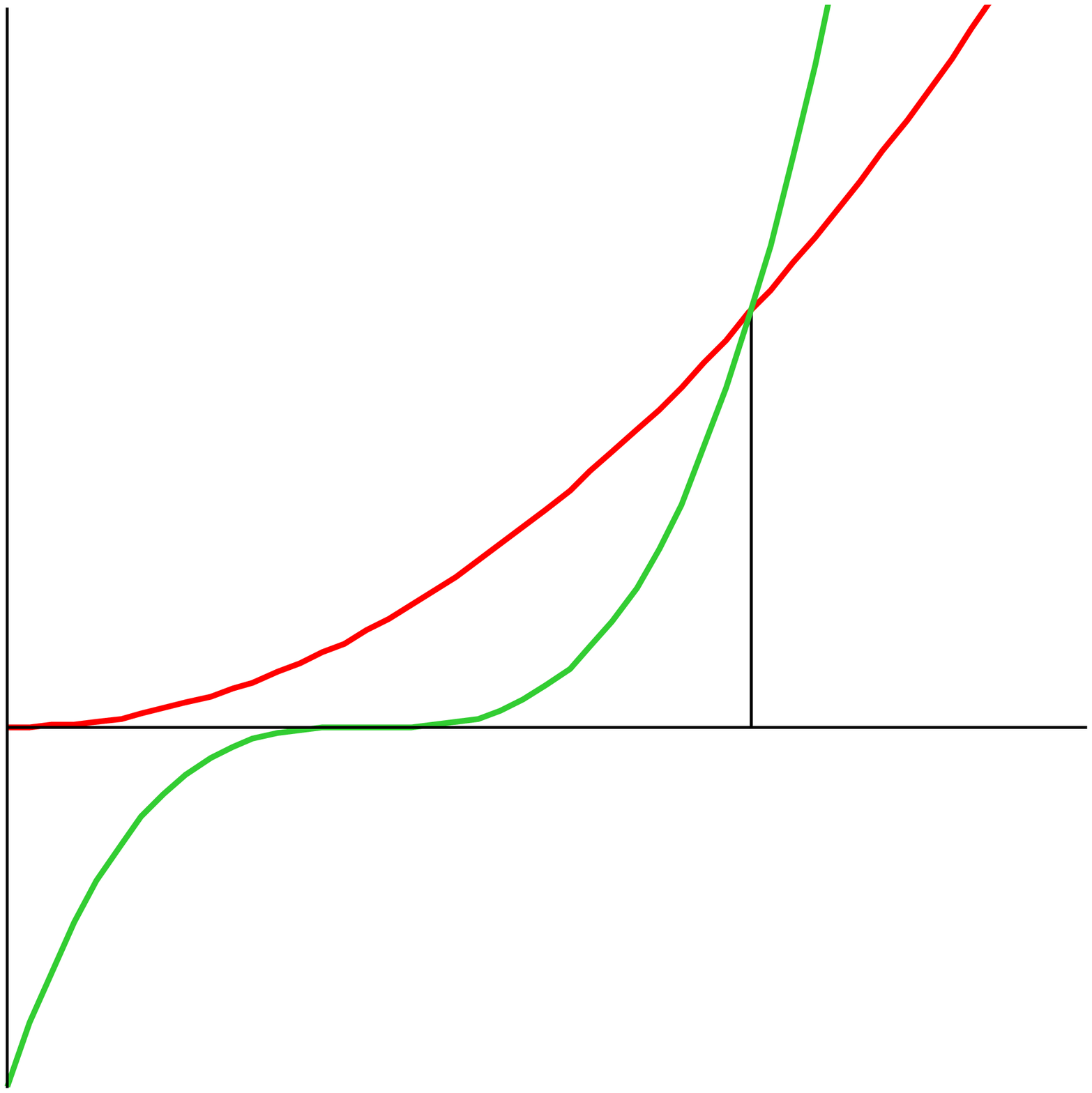}}}
\put(3.5,0){\rotatebox{0}{\includegraphics[scale=0.16]{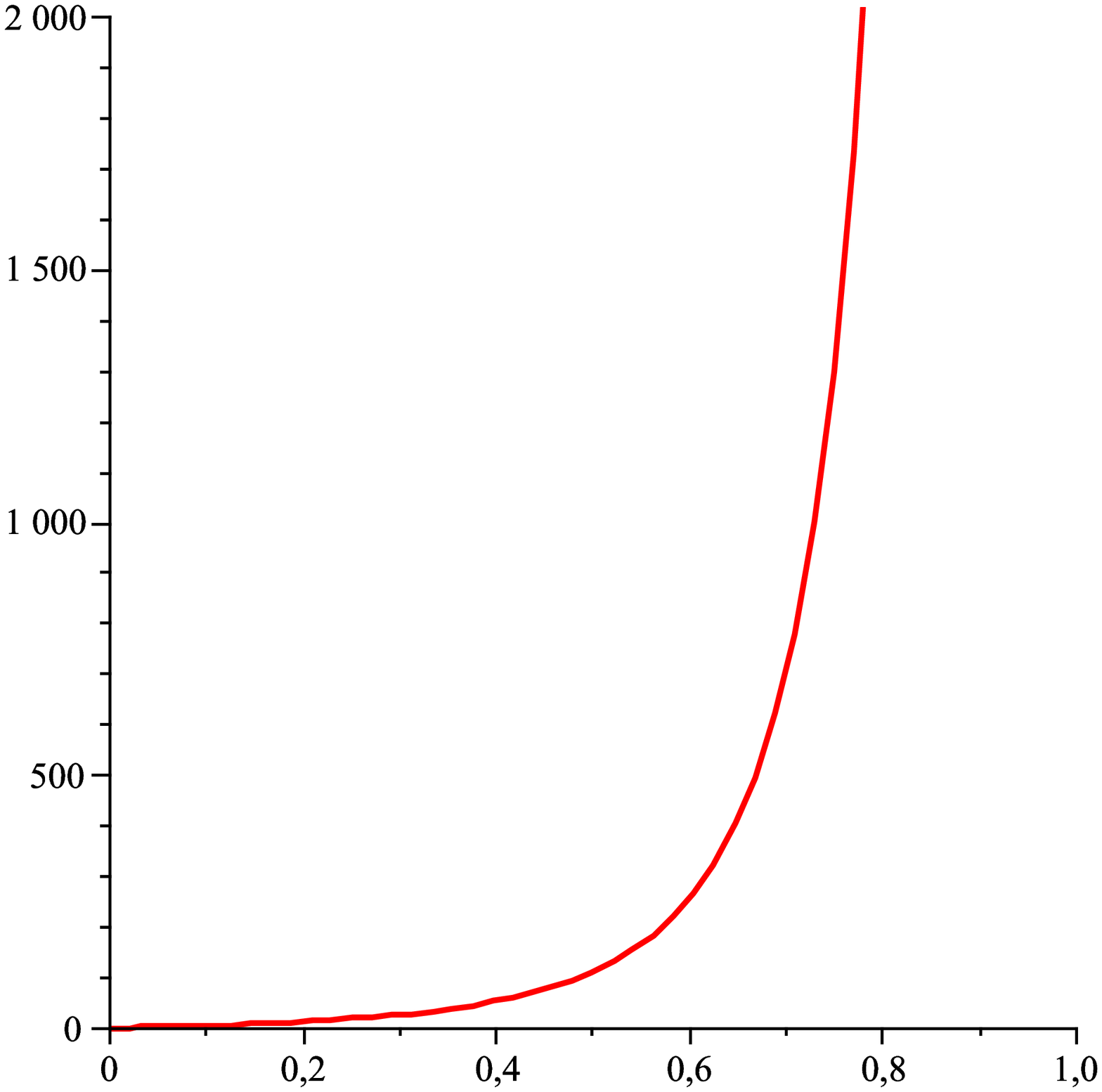}}}
\put(7,0){\rotatebox{0}{\includegraphics[scale=0.16]{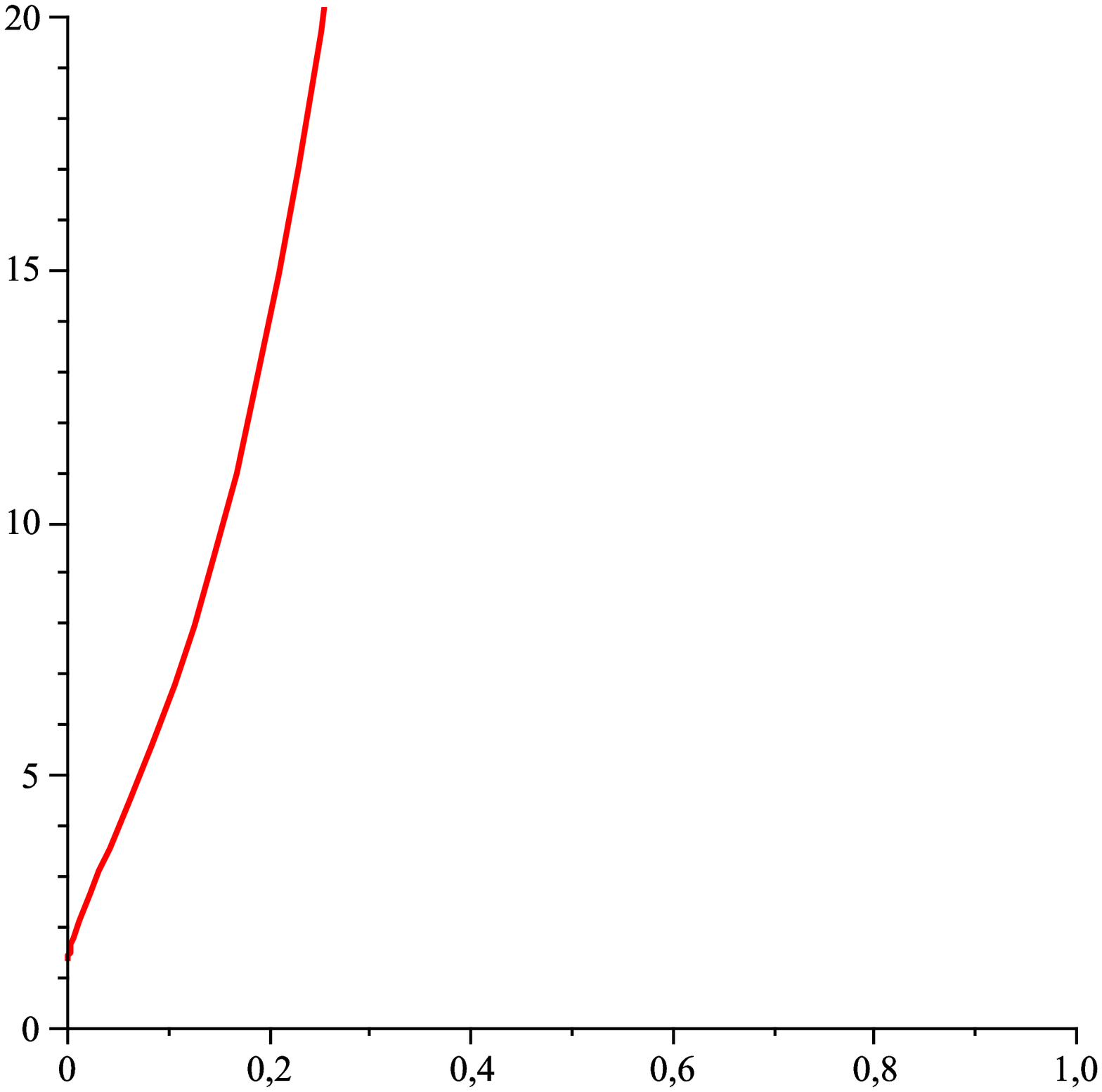}}}
\put(10.5,0){\rotatebox{0}{\includegraphics[scale=0.16]{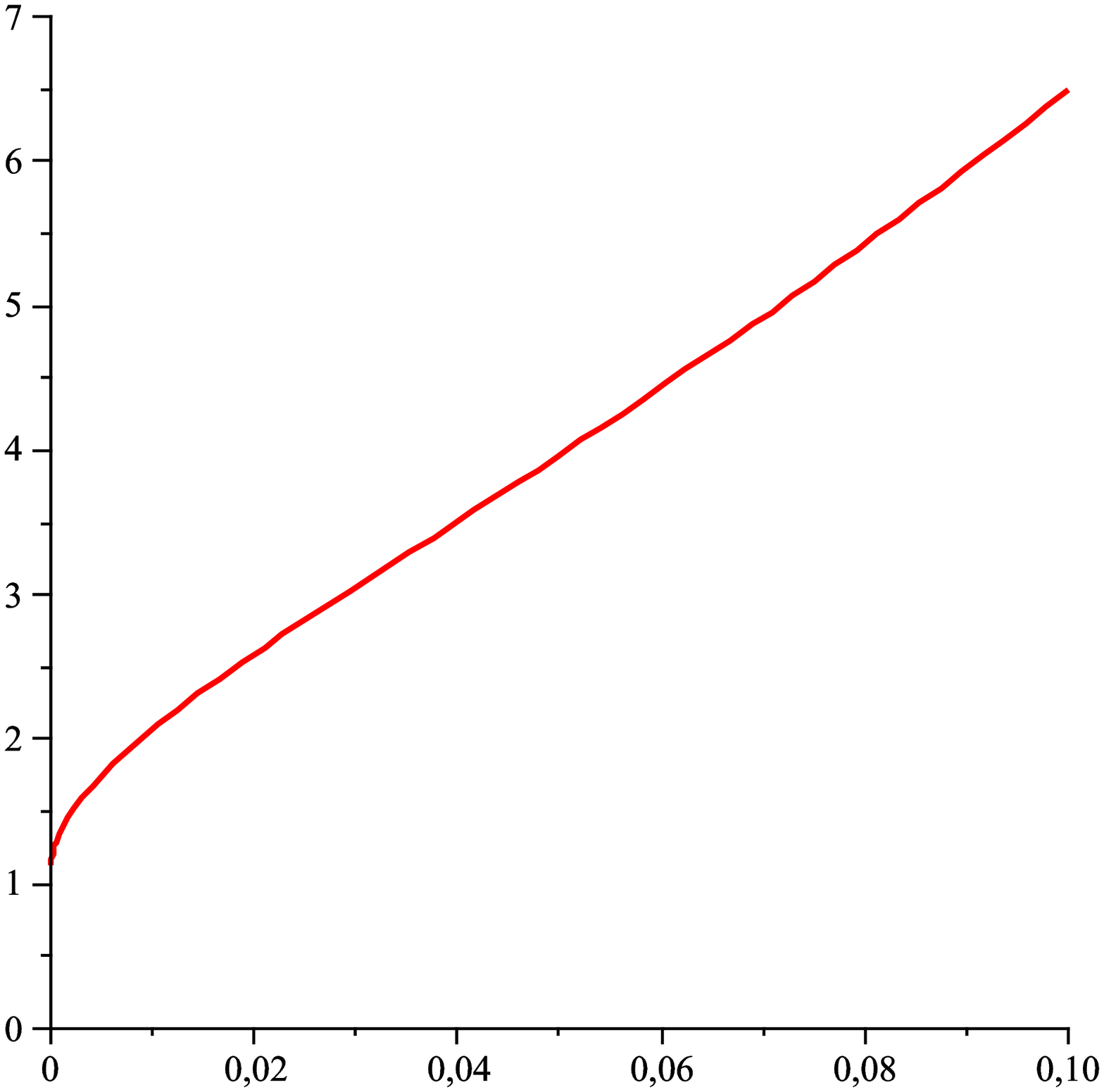}}}
\put(0.4,0.7){{\small $\frac{1}{1-b}$}}
\put(2.5,0.9){{\small $c$}}
\put(-0.65,2.8){{\small $y$}}
\put(1.3,0.8){{\small $c_{crit}(b)$}}
\put(2.2,2.8){{\small $y=27bc^2$}}
\put(-0.2,3){{\small $y=\left[c(1-b)-1\right]^3$}}
\end{picture}
\end{center}
\caption{On the left: the definition of the function $c_{crit}(b)$. 
For each $b<1$, 
the functions $y=\left[c(1-b)-1\right]^3$ (in green) 
and  $y=27bc^2$ (in red) intersect for $c=c_{crit}(b)$. 
On the right, the numerical plots (in red) of the function $c=c_{crit}(b)$, showing the quick increasing of $c_{crit}(b)$ with $b$.
} \label{figccrit}
\end{figure}
We need the following technical result.
\begin{lemma}\label{lemcrit}
The function $Q(S)=\frac{(1-S)(b+S)(1+cS)}{S}$ is decreasing over $[0,1]$
if and only if $$\left[c(1-b)-1\right]^3\leq27bc^2.$$
This condition is equivalent to either
$b\geq 1$ or $b<1$ and $c\leq c_{crit}(b)$,
where $c_{crit}(b)$ is the positive zero of $\left[c(1-b)-1\right]^3=27bc^2$. 
\end{lemma}
{\em Proof}
Since
$$Q'(S)=-\frac{2cS^3+\left(1+c(b-1)\right)S^2+b}{S^2},\quad Q''(S)=-\frac{2\left(b-cS^3\right)}{S^2},$$
the function $Q(S)$ has an inflexion point for $S=(b/c)^\frac13$. The function $Q(S)$ is nonincreasing
over $[0,1]$ if and only if its derivative at the inflexion point is nonpositive, 
that is, $P'\left((b/c)^\frac13\right)\leq0$.
Straightforward computations show that this condition is equivalent to $\left[c(1-b)-1\right]^3\leq27bc^2$. 
If $b\geq1$ then the first term of the inequality is negative and hence the inequality if satisfied for all $c\geq 0$. 
If $b<1$, then the inequality is satisfied if and only if $c\leq c_{crit}(b)$, see Fig. \ref{figccrit}.
The expression of $c_{crit}(b)$ can be obtained by Cardan formulas. Notice that $c_{crit}(0)=1$
and $c_{crit}(b)$ is quickly increasing with $b$.
\medskip
 
\begin{theorem}\label{ourthm2} 
Assume that 
\begin{equation}\label{1ppi}
\lambda_1<\lambda_2\leq\cdots\leq\lambda_N,\qquad \lambda_1<1,
\end{equation} 
\begin{equation}\label{ccrit1}
\mbox{either }b_1\geq 1\mbox{ or for each }i\geq 1\mbox{ satisfying }\lambda_i<1,~c_i\leq c_{crit}(b_1).
\end{equation}
Then the equilibrium $E_{1}^*$ is globally asymptotically 
stable for (\ref{eqMonodLQ}) with respect to the interior of the positive cone.
\end{theorem}
{\em Proof}
Let us prove that (\ref{h11}), (\ref{h31}) and (\ref{h2graphic}) hold.
The Monod function $f_1(S)$ is increasing. Hence, (\ref{h11}) holds. The function 
$P_1(S)$ is 
$$P_1(S)=\frac{(1-S)(b_1+S)(1+c_1S)}{S}.$$
By Lemma \ref{lemcrit}, it is decreasing if and only if either $b_1\geq 1$ or $b_1<1$ and $c_1\leq c_{crit}(b_1)$.
Hence, (\ref{h31}) holds.
For each $i\geq 2$, the function $g_i(S)$ defined by (\ref{gi}) is
$$g_i(S)=\frac{f_i(S)}{f_1(S)}\frac{1-S}{p_i(S)}
=\frac{Y_i}{a_i}
\frac{a_i-D_i}{a_1-D_1}
\frac{S-\lambda_i}{S-\lambda_1}Q_i(S) 
$$
where $Q_i(S)=\frac{(1-S)(b_1+S)(1+c_iS)}{S}$.
Assume that (\ref{ccrit1}) holds.
By Lemma \ref{lemcrit}, the function $Q_i(S)$ is decreasing. Therefore,
$$
\min_{0<S\leq\lambda_1}Q_i(S)=Q_i(\lambda_1)>Q_i(\lambda_i)=\max_{\lambda_i\leq S<1}Q_i(S).
$$
Since $\lambda_1<\lambda_i$, the function $S\mapsto \frac{S-\lambda_i}{S-\lambda_1}$ is increasing. Therefore,
$$\min_{0<S\leq\lambda_1}\frac{S-\lambda_i}{S-\lambda_1}=\frac{\lambda_i}{\lambda_1}>1>
\frac{1-\lambda_i}{1-\lambda_1}=\max_{\lambda_i\leq S<1}\frac{S-\lambda_i}{S-\lambda_1}.$$
Thus,
$$\min_{0<S<\lambda_1}g_i(S)\geq
\frac{Y_i}{a_i}\frac{a_i-D_i}{a_1-D_1}\min_{0<S<\lambda_1}\frac{S-\lambda_i}{S-\lambda_1}\min_{0<S<\lambda_1}Q_i(S)
>\frac{Y_i}{a_i}\frac{a_i-D_i}{a_1-D_1}Q_i(\lambda_1),$$
and 
$$\max_{\lambda_i<S<1}g_i(S)\leq
\frac{Y_i}{a_i}\frac{a_i-D_i}{a_1-D_1}\max_{\lambda_i<S<1}\frac{S-\lambda_i}{S-\lambda_1}\max_{\lambda_i<S<1}Q_i(S)<
\frac{Y_i}{a_i}\frac{a_i-D_i}{a_1-D_1}Q(\lambda_i)
.$$
Hence (\ref{h2graphic}) holds. The result follows from Theorem \ref{ourthm1}.
\medskip

\begin{figure}[ht]
\setlength{\unitlength}{1.0cm}
\begin{center}
\begin{picture}(9,3.8)(0,0)
\put(0,0){\rotatebox{0}{\includegraphics[scale=0.2]{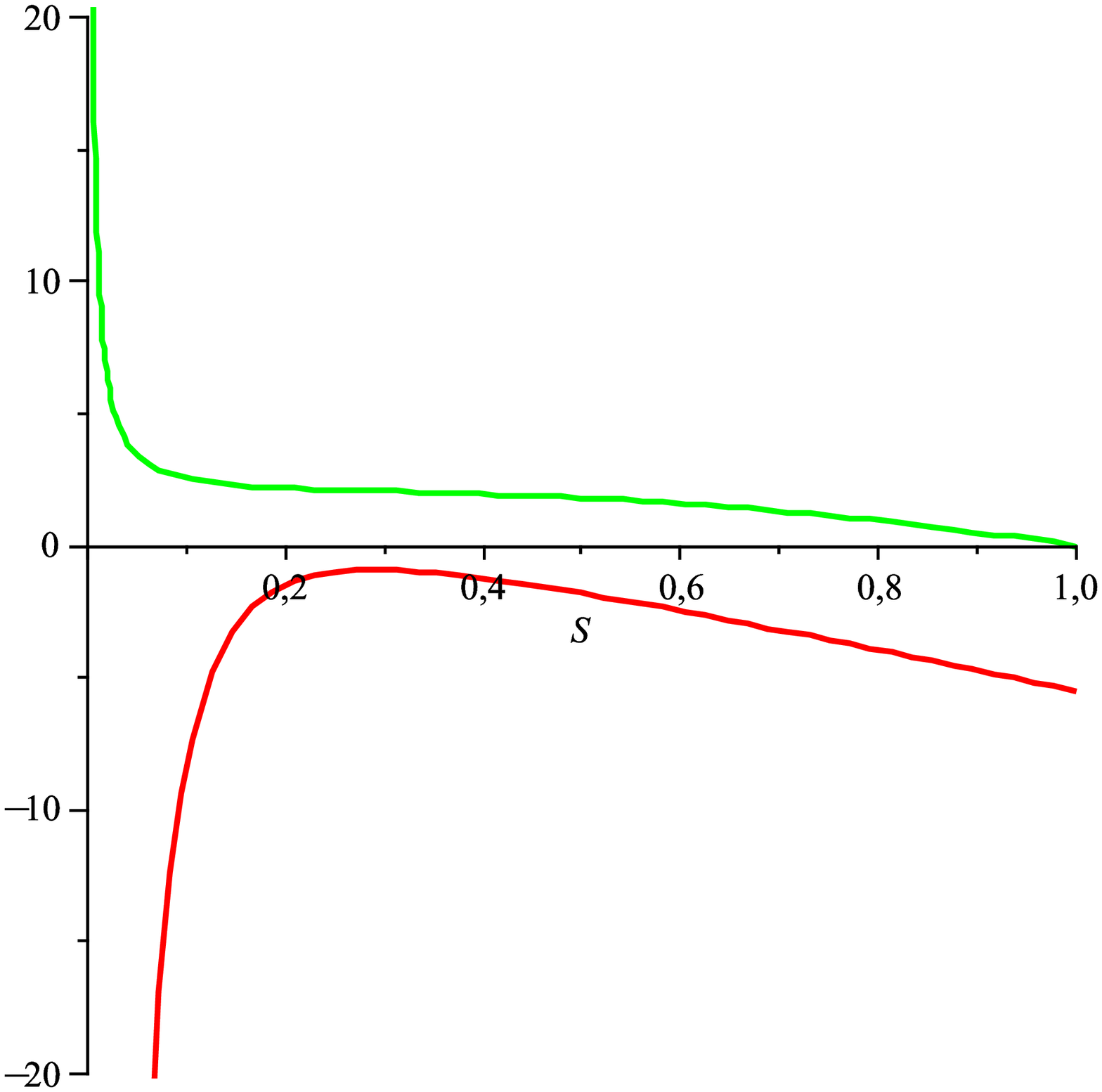}}}
\put(4.5,0){\rotatebox{0}{\includegraphics[scale=0.2]{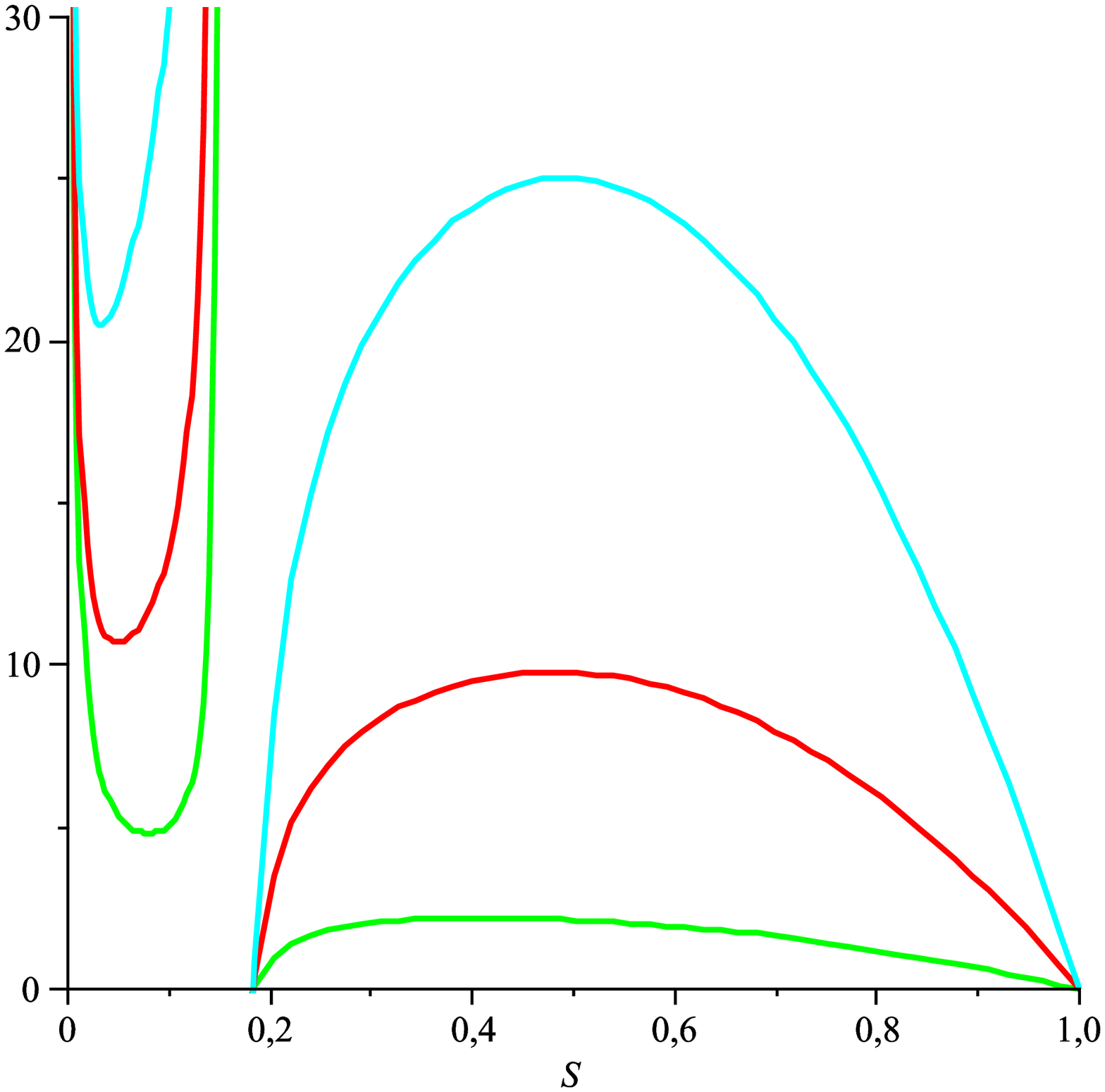}}}
\end{picture}
\end{center}
\caption{The graphical depiction of conditions (\ref{h31}) and (\ref{h2graphic}) in the proof of Theorem  \ref{ourthm2}.
The parameter values are
$c_1 = 4$, $b_1 = 0.1$, $a_1 = 1$,
$b_2 = 0.15$, $a_2 = 1$
$D_1 = 0.6$ and $D_2 = 0.55$. Hence
$\lambda_1=0.15$ and
$\lambda_2\simeq0.18$. 
On the left, the function $P_1(S)$ (in green) and its derivative (in red) showing that $P_1$
is decreasing and so (\ref{h21}) is satisfied.
On the right, the function $g_2(S)$ for 
$c_2 = 5$ (in red),
$c_2 = 30$ (in green) and
$c_2 = 80$ (in cyan). The condition (\ref{h2graphic}) is satisfied for  $c_2 = 5<c_{crit}(0.11)$ and
$c_2 = 30>c_{crit}(0.1)$. It is not satisfied for $c_2 = 80$. Here $c_{crit}(0.1)\simeq6.5$, see Fig. \ref{figccrit}}  \label{figF}
\end{figure}

Theorem  \ref{ourthm2} extends Corollary \ref{thmhsu} 
which corresponds to the case where the yields are constant. 
Indeed, (\ref{1ppi}) is the same as (\ref{hcm1}) and, for constant yields, $c_i=0$, 
so that the conditions (\ref{ccrit1}) in Theorem \ref{ourthm2} are satisfied. 
Notice that (\ref{ccrit1}) is a sufficient and not necessary condition for the existence 
of a gap between the minimum of $g_i(S)$ over $(0,\lambda_1)$ and its maximum over
$[\lambda_i,1]$. For instance, 
for the parameter values given in the caption of Fig. \ref{figF},
if $c_2=30>6.5\simeq c_{crit}(0.1)$, there exists such a gap, see Fig. \ref{figF}. 
Therefore, Theorem \ref{ourthm1} applies and predict
that the equilibrium is GAS, even if Theorem \ref{ourthm2} does not apply, 
since $b_1=0.1<1$ and $c_2=30>6.5\simeq c_{crit}(0.1)$. 
However, for $c_2=80>6.5\simeq c_{crit}(0.1)$, there is no gap, see Fig. \ref{figF}. 
Therefore, neither Theorem \ref{ourthm1} nor Theorem \ref{ourthm2} can be used.
This model shows that in each particular example it is very easy to depict graphically 
the conditions (\ref{h11}), (\ref{h31}) and (\ref{h21})
of Theorem \ref{ourthm1} and to see if Theorem \ref{ourthm2} can be applied.

\section{Discussion}\label{conclusions}

We briefly survey some CEP results for (\ref{eqsxi}). 
In the Monod case \cite{monod} when the growth functions are of form (\ref{monod}),
and assuming equal removal rates for $S$ and all species, 
i.e. $D_i=D$ for $i=1\cdots N$, 
Hsu, Hubbell and Waltman \cite{hhw} proved the following CEP: 
every solution of (\ref{eqsxi}) with positive initial condition satisfies
$$\lim_{t\to\infty}S(t)=\lambda_1,
\quad \lim_{t\to\infty}x_1(t)=Y_1(S^0-\lambda_1),
\quad
\lim_{t\to\infty}x_i(t)=0,~i\geq 2,
$$
under the additional assumption $0<\lambda_1<S^0$ and 
$\lambda_1<\lambda_i$ for $i=2\cdots N$.
The predictions in \cite{hhw} were tested in the laboratory 
by the experiments of Hansen and Hubbell \cite{hh}. 
Similar experiments could be performed
to test the predictions using microorganisms known to have variable yields. 
See \cite{ALLR,apw,CT,CWT,PW} for examples of such microorganisms.
 
Hsu \cite{hsu} used the Lyapunov function (\ref{HSU})  
to give a simple and elegant proof of the result in \cite{hhw} 
for the case of different removal rates $D_i$ (see Corollary \ref{thmhsu}).
Wolkowicz and Lu \cite{wl}, used the Lyapunov function  (\ref{LWL}) 
and extended the results of \cite{hsu} by allowing 
more general growth functions. 
They identified a large class of growth functions, where
the constant $c_i$ in (\ref{LWL}) can always be found. 
Following Smith and Waltman \cite{chem}, I think that,
despite the fact the  $c_i$ cannot be found for all growth functions, 
the work of Wolkowicz and Lu \cite{wl} represents a major step
in the extension of the result of Hsu \cite{hsu} to general growth functions. 
In the constant yield case, the CEP has been also proved under a variety of hypotheses by
Armstrong and McGehee \cite{amg},
Butler and Wolkowicz \cite{bw}, Wolkowicz and Xia \cite{wx}
and Li \cite{li}. The hypotheses used in the papers \cite{amg,bw,hsu,hhw,li,wl,wx} 
are summarized in Table 1 of \cite{LLS}. Lyapunov techniques in the chemostat 
were also used in \cite{LM,RH}.

The variable yield case was considered, for $n=1,2$ by Pilyugin and Waltman \cite{PW}, 
with a particular interest to linear and quadratic yields, and by 
Huang, Zhu and Chang \cite{HZC}. 
The general model (\ref{eqsqi}) for $N$ species, was considered by 
Arino, Pilyugin and Wolkowicz \cite{apw}.
As noticed by these authors (see \cite{apw}, Section 3), in the case of constant yields (\ref{eqsxi}), 
including the yield terms $Y_i$  in the substrate equation, as in (\ref{eqsxi}), is mathematically
equivalent to including the reciprocal in the microorganism equation instead.
Indeed, (\ref{eqsxi}) can be written
\begin{equation*}
\begin{array}{lcll} 
 \displaystyle S'& =&D[S^0-S]-\displaystyle
\sum_{i=1}^N{p_i(S)}x_i,&\\[3mm]
 \displaystyle x'_i&=& [Y_ip_i(S) - D_i]x_i,& i=1\cdots N,
\end{array}
\end{equation*}
where $p_i(S)=\frac{q_i(S)}{Y_i}$. 
Since $Y_i$ are constant, the uptake terms $p_i$ and
growth terms $q_i$ have the same monotonicity properties.  
Formally, the model (\ref{eqsqi}-\ref{yields}) with variable yields $y_i(S)$ can be written  
\begin{equation}\label{Ey1}
\begin{array}{lcll} 
 \displaystyle S'& =&D[S^0-S]-\displaystyle
\sum_{i=1}^N\frac{q_i(S)}{y_i(S)}x_i,&\\[3mm]
 \displaystyle x'_i&=& [q_i(S) - D_i]x_i,& i=1\cdots N,
\end{array}
\end{equation}
where $q_i(S)$ are the growth functions, or equivalently,
\begin{equation}\label{Ey2}
\begin{array}{lcll} 
 \displaystyle S'& =&D[S^0-S]-\displaystyle
\sum_{i=1}^Np_i(S)x_i,&\\[3mm]
 \displaystyle x'_i&=& [y_i(S)p_i(S) - D_i]x_i,& i=1\cdots N,
\end{array}
\end{equation}
where $p_i(S)=\frac{q_i(S)}{y_i(S)}$ are the uptake functions.
One of the important differences in the case that the yields are not constant 
is that the variable yield terms can lead to uptake and
growth functions that have now different monotonicity properties. 
Moreover, in the case of constant yields, the yields terms $Y_i$ can be eliminated
in (\ref{eqsxi}) simply by passing to the variables $u_i=Y_ix_i$. We obtain 
\begin{equation*}
\begin{array}{lcll} 
 \displaystyle S'& =&D[S^0-S]-\displaystyle
\sum_{i=1}^N{q_i(S)}u_i,&\\[3mm]
 \displaystyle u'_i&=& [q_i(S) - D_i]u_i,& i=1\cdots N.
\end{array}
\end{equation*}
This change of variables means that we have changed the units in which 
the microorganisms were evaluated. 
There is no such trick to eliminate the yields terms in (\ref{Ey1}) of (\ref{Ey2}).
Therefore, careful attention to the interpretation of the yield terms 
resulting in their correct placements in the equations is necessary. 
For details and complements, the reader is referred to \cite{apw}, Section 3.

In the variable yield case, the CEP has been proved for (\ref{eqsqi}), under some technical 
conditions on the function $p_i$ and $q_i$, 
by Sari \cite{sari1} and Sari and Mazenc \cite{sari2} (see Corollary \ref{ourthms} in Section \ref{wolkowicz}
and Corollary \ref{ourthmsm}  in Section \ref{hsu}). 
It was also shown in \cite{sari2} 
how Corollary \ref{ourthmsm} can be fruitfully used to analyze the stability properties of systems whose yield 
functions depend on the variable $S$. For instance, 
the CEP holds (see Corollary 5 in \cite{sari2})
for the Monod model with constant yields replaced by either linear or 
quadratic functions of $S$, and under certain additional technical assumptions. 
Another application is given by the model of Pilyugin and Waltman \cite{PW} 
which was used to demonstrate that a periodic orbit was possible in the case of variable yield model. 
In this model, with two species, where one yield is constant and the other is cubic in $S$, 
it is shown in \cite{sari2} that for some values of the parameters
the CEP holds (see Corollary 6 in \cite{sari2}). The problem of the existence of limit cycles in chemostat 
equations is not always well understood \cite{sari}.
In the case of constant yields, numerical simulations 
of model (\ref{eqsxi}) have only displayed 
competitive exclusion. Our results concern also
the case of variable yields, for which it is known \cite{apw,HZC,PW} that more exotic dynamical behaviors, 
including limit cycles and chaos, are possible. Thus in the case of variable yields, it is of great importance 
to have criteria ensuring the global convergence to an equilibrium with at most one surviving species.
The reader interested in biological motivations for the dependence of the yields 
on the substrate, may consult \cite{apw,PW} and the references therein.

The exclusion of periodic orbits in system (\ref{D=1}) was obtained 
by Fiedler and Hsu (see Theorem 1.1 in \cite{FH})
under the following conditions: for all $1\leq i\neq j\leq N$,
and $0<S<1$ 
\begin{equation}\label{FH1}
(S-\lambda_i)f_i(S)>0,\quad\mbox{ for }S\neq \lambda_i,
\end{equation}
\begin{equation}\label{FH2}
f_i(S)<1+f_j(S)+(1-S)p'_j(S)/p_j(S).
\end{equation}
Even if the result in \cite{FH} does not show the convergence to an equilibrium, 
it is interesting to compare the constraints on the functions $f_i$ and $p_i$ 
of \cite{FH} with our constraints.
Actually, (\ref{FH1}) is stronger than (\ref{h11}), 
since our assumption requires (\ref{FH1}) only for $f_1$ 
and allows the $f_i$, for $i\neq 1$, 
to have other zeros than $\lambda_i$ in $]0,1[$. 
Let us compare the constraints on  $f_i$ and $p_i$ imposed by the inequalities (\ref{FH2})
to the constraints imposed by hypothesis (\ref{h21}): notice that (\ref{FH2}) 
is a set of $\frac{N(N-1)}{2}$ conditions, while (\ref{h21}) is a set of at most $N-1$ conditions. 
Moreover, the constants $\alpha_i$ in the conditions (\ref{h21}) give more flexibility to 
these conditions. For instance, (\ref{h21}) 
are satisfied by arbitrary Monod growth functions and also by a large class of growth 
functions as it was shown in \cite{sari2,wl}. 
On the other hand, (\ref{FH2}) are not satisfied by arbitrary
Monod functions, see formulas (6.10) and (6.11) in \cite{FH}. 
Hence the result in \cite{FH} does not recover the CEP, even in the 
classical and well established case of Monod functions and equal removal rates \cite{amg}. 
However our theorem recovers a lot of
results of the existing literature. 
The beautiful geometrical techniques in \cite{FH}, 
whose purpose is to extend Bendixon-Dulac criterion to higher dimension are innovative 
and it seems likely that \cite{FH} is destined to be cited more often 
for these innovative techniques rather than for the specific results 
that are proved for the CEP. 

For the purpose of comparison between our result and the 
result of Fiedler and Hsu \cite{FH}, we just mentioned  
two caveats on Theorem 1.1 in \cite{FH}: 
first, this theorem does not recover many of the biologically 
interesting classical examples where the CEP is known to hold, 
and second, it does not prove the convergence to an equilibrium. 
These caveats were already mentioned in \cite{FH}, Section 6.
Another caveat must be signaled. 
Fiedler and Hsu claimed (see \cite{FH}, Section 6) that,
in the case $N=1$ of a single 
species, condition (\ref{FH2}) holds trivially and there is no periodic orbit for
system (\ref{eqsx}). It should be noticed that condition (\ref{FH1}) 
is not sufficient to exclude periodic orbits. 
Of course, if $p'(S)>0$, then the Bendixon-Dulac criterion can be applied to exclude
periodic orbits (see Section \ref{N=1}). This assumption on the monotonicity of $p$
is not explicitly stated in Theorem 1.1 in \cite{FH}. 
Moreover, the condition $p'(S)>0$ would not be satisfactory 
from the biological point of view. 
Indeed, a variable yield term $y(S)=q(S)/p(S)$ 
can lead to nonmonotone uptake term $p(S)$ even if the
growth term $q(S)$ is monotone (see Section \ref{N=1}). 

Fiedler and Hsu, see Section 6 in \cite{FH}, 
claimed that the construction of Lyapunov functions in  \cite{amg,hsu,li,wl,wx} 
strictly depends on the proportionality (\ref{proportionality}) required in equations (\ref{eqsxi}).
We show how the Lyapunov function used by Hsu himself \cite{hsu} for the Monod case,  
more than thirty years ago, can be extended to the case of (\ref{eqfpi}), 
where growth rates are not required to be proportional to food uptake
(see Theorem \ref{ourthm}, in Section \ref{hsu}).
For that reason, the direct proof of Theorem \ref{ourthm}, using 
the extension of the Lyapunov function of Hsu \cite{hsu}, seems to be interesting in itself.  
Thus, I decided to give  Theorem \ref{ourthm} and its 
direct proof, despite the fact that this theorem is a corollary 
of Theorem \ref{ourthm1} (see Proposition \ref{prop1} in Section \ref{hsu}).

We list some references to the existing literature which inspired our approach.
The Lyapunov function (\ref{lyapunov}) used in the proof of Theorem \ref{ourthm} 
was introduced in \cite{sari2} as an extension of the 
Lyapunov function (\ref{HSU}) 
that Hsu used in \cite{hsu} in the Monod case (Theorem 3.3 in \cite{hsu}).
In the case of one species, this Lyapunov function is equal to the   
function (\ref{LPW}) used by Pilyugin and Waltman (Lemma 2.3 in \cite{PW}), 
as shown in Section \ref{N=1}.
It is also a multiple of the Lyapunov function that Ballyk, Lu, Wolkowicz and Xia used 
in \cite{wx}, page 1039 or \cite{wbl}, Section 3.3 (see Section 3.2 in \cite{sari2}).
The Lyapunov function (\ref{lyapunov1}) used in the proof of Theorem \ref{ourthm1} 
was introduced in \cite{sari1} as an extension of the Lyapunov function (\ref{LWL}) 
that Wolkowicz and Lu used in \cite{wl} in the constant yields case (Theorem 2.3 in \cite{wl}). 
In the case of one species, it is a multiple of the Lyapunov  
function (\ref{APW}) used by Arino, Pilyugin and Wolkowicz (Theorem 2.11 in \cite{apw}), 
as shown in Section \ref{N=1}.

In this work, we have analyzed a general model of the chemostat with several 
species competing for a substrate,
under the assumption that uptake rates and growth rates are not proportional.
Each species is characterized by its specific growth rate, its specific removal rate, 
and its variable yield.  
Our study reveals that the CEP holds for a large class of systems: 
the species with the smallest break-even concentration can be the winner of the competition
if some supplementary conditions, involving the uptake and growth functions are satisfied. 
Hence, even if the break-even concentration are depending only on the
growth rates and not on the yields functions, 
the issue of competition depends really on the yield functions. 
For instance, if one on the species exhibits a linear yield, and if the parameter
in the yield is enlarged, then the equilibrium, where only the winning species survives,
can be destabilized, and oscillatory coexistence of more than one species becomes possible.

\section*{Acknowledgments} 
The author gratefully acknowledge Alain Rapaport and J\'er\^ome Harmand for fruitful discussions.

\end{document}